\documentclass[review,hidelinks,onefignum,onetabnum]{siamart220329}

\usepackage{lipsum}
\usepackage{amsfonts}
\usepackage{graphicx}
\usepackage{epstopdf}
\usepackage{algorithmic}
\ifpdf
  \DeclareGraphicsExtensions{.eps,.pdf,.png,.jpg}
\else
  \DeclareGraphicsExtensions{.eps}
\fi

\newsiamremark{remark}{Remark}
\newsiamremark{hypothesis}{Hypothesis}
\crefname{hypothesis}{Hypothesis}{Hypotheses}
\newsiamthm{claim}{Claim}

\headers{Generalized Optimal AMG Convergence Theory}{Ali, Brannick, Kahl, Krzysik, Schroder, Southworth}

\title{Generalized Optimal AMG Convergence Theory for Nonsymmetric and Indefinite Problems\thanks{Submitted to the editors July 22, 2024.
\funding{This work was funded by NSF grants DMS-2110917 and DMS-2111219. B. S. S. was supported by the Laboratory Directed Research and Development program of Los Alamos National Laboratory under project number 20240261ER. LA-UR-24-26359.}}}

\author{Ahsan Ali\thanks{Department of Mathematics and Statistics, University of New Mexico, Albuquerque, NM 87131, USA (\email{ahsan@unm.edu}).}
\and James J. Brannick\thanks{Department of Mathematics, Penn State University, University Park, PA 16802, USA (\email{jjb23@psu.edu}).}
\and Karsten Kahl\thanks{School of Mathematics and Natural Sciences, Bergische Universit\"{a}t Wuppertal, Wuppertal 42119, Germany (\email{kkahl@uni-wuppertal.de}).}
\and Oliver A. Krzysik\thanks{Theoretical Division, Los Alamos National Laboratory, Los Alamos, NM 87545, USA (\email{okrzysik@lanl.gov}).}
\and Jacob B. Schroder\thanks{Department of Mathematics and Statistics, University of New Mexico, Albuquerque, NM 87131, USA (\email{jbschroder@unm.edu}).}
\and Ben S. Southworth\thanks{Theoretical Division, Los Alamos National Laboratory, Los Alamos, NM 87545, USA (\email{southworth@lanl.gov}). }}

\usepackage{amsopn}

\ifpdf
\hypersetup{
  pdftitle={Generalized Optimal AMG Convergence Theory for Nonsymmetric and Indefinite Problems},
  pdfauthor={Ahsan Ali, James J. Brannick, Karsten Kahl, Oliver A. Krzysik, Jacob B. Schroder, and Ben S. Southworth}
}
\fi

\usepackage{lipsum}
\usepackage{graphicx}
\usepackage{caption,subcaption}
\usepackage{amsmath}
\usepackage{amssymb,mathtools}
\usepackage{amsfonts}
\usepackage{float}
\usepackage{cleveref}
\usepackage{booktabs}
\usepackage{comment}
\usepackage{bm}
\usepackage{xcolor}
\usepackage{enumerate}
\renewcommand{\theenumi}{\roman{enumi}}%

\usepackage{url}

\DeclareFontFamily{U}{mathx}{}
\DeclareFontShape{U}{mathx}{m}{n}{<-> mathx10}{}
\DeclareSymbolFont{mathx}{U}{mathx}{m}{n}
\DeclareMathAccent{\widehat}{0}{mathx}{"70}
\DeclareMathAccent{\widecheck}{0}{mathx}{"71}

\nolinenumbers

\begin{document}

\maketitle
\begin{abstract}
Algebraic multigrid (AMG) is known to be an effective solver for many sparse symmetric positive definite (SPD) linear systems. For SPD systems, the convergence theory of AMG  is well-understood in terms of the $A$-norm, but in a nonsymmetric setting such an energy norm is non-existent. For this reason, convergence of AMG for nonsymmetric systems of equations remains an open area of research. A particular aspect missing from theory of nonsymmetric and indefinite AMG is the incorporation of general relaxation schemes. In the SPD setting, the classical form of optimal AMG interpolation provides a useful insight in determining the best possible two grid convergence rate of a method based on an arbitrary symmetrized relaxation scheme. In this work, we discuss a generalization of the optimal AMG convergence theory targeting nonsymmetric problems, using a certain matrix-induced orthogonality of the left and right eigenvectors of a generalized eigenvalue problem relating the system matrix and relaxation operator.  We show that using this generalization of the optimal convergence theory, one can obtain a measure of the spectral radius of the two grid error transfer operator that is mathematically equivalent to the derivation in the SPD setting for optimal interpolation, which instead uses norms. In addition, this generalization of the optimal AMG convergence theory can be further extended for symmetric indefinite problems, such as those arising from saddle point systems so that one can obtain a precise convergence rate of the resulting two-grid method based on optimal interpolation.  We provide supporting numerical examples of the convergence theory for  nonsymmetric advection-diffusion problems, the two-dimensional Dirac equation motivated by $\gamma_5$-symmetry, and the mixed Darcy flow problem corresponding to a saddle point system.
\end{abstract}



\section{Introduction}\label{sec:intro} Algebraic Multigrid (AMG) \cite{brandt1986algebraic,brandt1984algebraic,JWRuge_KStuben_1987a} is a powerful and efficient iterative solver for solving large-scale linear systems arising from various applications, including computational fluid dynamics, structural mechanics, and other scientific simulations. When typical discretization approaches such as finite difference, finite volume, or finite element are applied in the numerical approximation to partial differential equations, the resulting approximation can swiftly grow to several million, billions, or even more, numbers of unknowns. AMG methods are often used to solve such large sparse linear system of equations\begin{equation}\label{eq:main}
A \bf{x} = \bf{b},
\end{equation}
where $A \in \mathbb{C}^{n\times n}$
 is a sparse nonsingular matrix with $\bf{x},\bf{b}$ $ \in \mathbb{C}^{n}$. AMG is one of the most widely used linear solvers due to its potential for $\mathcal{O}(n)$ computational cost when solving a linear system with $n$ number of degrees of freedoms.

Convergence of AMG is well-motivated when $A$ is symmetric positive definite (SPD) and induces an energy norm, which frequently arises from the discretization of elliptic or certain parabolic partial differential equations (PDEs) \cite{falgout2004generalizing,falgout2005two,maclachlan2014theoretical,xu2017algebraic,notay2015algebraic,vassilevski2008multilevel,brannick2018optimal}. When $A$ is nonsymmetric, the AMG framework has been extended to nonsymmetric matrices using a number of different techniques~\cite{sala2008new, brezina2010towards,brannick2014bootstrap,lottes2017towards,manteuffel2017root,manteuffel2018nonsymmetric,ali2024constrained,wiesner2014multigrid}. The recent work \cite{manteuffel2019convergence} made progress towards nonsymmetric AMG theory utilizing a generalized ``fractional'' approximation property to demonstrate convergence. Two-level theory based on approximating so-called ideal transfer operators was also developed in \cite{manteuffel2018nonsymmetric}. In \cite{lottes2017towards}, the author established nonsymmetric AMG theory by forming the absolute value matrix $|A|$.  The work~\cite{southworth2023compatible} developed compatible transfer operators in nonsymmetric AMG with standard matrix-induced norms. Thus, while substantial progress has been made in the development of nonsymmetric AMG solvers, there is still a need for reliable and predictive theoretical motivation to investigate the robustness of existing, as well as new methods. In particular, recent progress on two-level theory for nonsymmetric AMG \cite{manteuffel2018nonsymmetric,manteuffel2019convergence,lottes2017towards,southworth2023compatible} does not incorporate the effect of general relaxation operators or the necessary relationship between multigrid transfer operators and the relaxation scheme, whereas both of these issues are accounted for by our theory.
Overall, the incomplete theoretical framework in the nonsymmetric setting contributes to lack of robustness for nonsymmetric solvers, when compared to the SPD case, and thus motivates this work. Robustness in this context refers to the solver's ability to achieve consistent and reliable convergence rates across a wide range of nonsymmetric problems, and traditionally a firm theoretical understanding of multigrid has aided in the development of such solvers.

The efficacy of AMG methods as solvers depends on the complementary processes of relaxation and coarse-grid correction, both of which aim to reduce the error in a solution iteratively. Here, relaxation refers to a simple iterative technique (e.g. Gauss-Seidel) that efficiently reduces some part of the error, e.g., in the classic SPD diffusion case, high-frequency error components associated with rapidly oscillating modes. Meanwhile, coarse-grid correction tackles the remaining error, e.g., in the classic case, low-frequency components left over by relaxation. Together, these two components ensure that all aspects of the error are addressed and enable AMG to act as an efficient iterative solver.

To set up a coarse-grid correction, we first partition $A$, which represents the system matrix of the problem to be solved, into two groups of unknowns: fine points ($F$) and coarse points ($C$). This partition, known as the $C/F-$splitting, divides the $n$ unknowns in $A$ such that there are $n_f$ fine nodes in $F$ and $n_c$ coarse nodes in $C$, with $C\cap F = \emptyset$ and $C\cup F = \{1,2,\dots,n\}$. This process is foundational in AMG, as it allows us to separate the matrix into components that are treated differently in the solution process.

Reordering the matrix $A$ with $F$-points followed by $C$-points yields
\begin{equation}\label{CFspl}
A = \left[\begin{array}{cc}
A_{f\!f}   & A_{f\!c} \\
A_{c\!f} & A_{cc} \\
\end{array}\right],
\end{equation}
where for example, $A_{f\!f}$ corresponds to entries in $A$ between two $F$-points. This structure is crucial for constructing the transfer operators, which map information between the fine and coarse levels, allowing error reduction across different scales.

The error reduction process in AMG is encapsulated in the two-grid error transfer operator, which combines relaxation and coarse-grid correction to form a single operator, denoted by $E_{TG}$, that describes how the error evolves from one iteration to the next \begin{equation}\label{ETG}
    E_{TG}=(I-M^{-*}A)(I-PA_{c}^{-1}RA)(I-M^{-1}A).
\end{equation}
Here, $R:\mathbb{C}^{n}\mapsto\mathbb{C}^{n_{c}} (n_{c}< n)$ denotes the restriction operator (maps the residual from the fine level to the coarse level), $P:\mathbb{C}^{n_{c}}\mapsto\mathbb{C}^{n}$ denotes the interpolation matrix (maps error correction from the coarse level back to the fine level), $A_{c}=RAP$ is the coarse-grid operator, $M$ denotes the relaxation, $M^{*}$ denotes the conjugate transpose of $M$, and $(I-PA_{c}^{-1}RA)$ denotes the error transfer operator of just coarse-grid correction. We use the notation $M^{-*}$ to denote the inverse conjugate transpose of $M$, which is $(M^{*})^{-1} = (M^{-1})^{*} \coloneqq M^{-*} $. In any specific case, if $M$ is real, then $M^{*}$ is just the transpose of $M$. The action of $E_{TG}$ defines how the error is reduced through each iteration. Ideally, the norm of $E_{TG}$ should be less than one in some suitable norm, ensuring that the error diminishes across iterations. In fact, two-level multigrid theory largely concerns itself with the analysis of and attempts to bound $E_{TG}$. If the AMG process is done recursively, i.e., AMG is applied to $RAP$, then a multilevel hierarchy is achieved. 

The error propagator \eqref{ETG} suggests that classical AMG was initially intended for SPD problems. In this instance, the coarse grid correction $I-PA_{c}^{-1}RA$ becomes an $A$-orthogonal projector when $R = P^{*}$. The two-grid approach is convergent in the $A$-norm if it is combined with a convergent global relaxation scheme $M$. 
Here, global relaxation refers to a relaxation scheme that updates all degrees of freedom and includes methods like standard Jacobi or Gauss-Seidel.
While AMG works effectively on many SPD problems \cite{brandt1983algebraic,briggs2000multigrid,xu2017algebraic,JWRuge_KStuben_1987a}, the nonsymmetric case is far more challenging for comparable theory. As an example, $I-PA_{c}^{-1}RA$ is no longer an orthogonal projection and will cause an increase in error in any standard inner-product. Furthermore, the construction of convergent relaxation methods becomes more complicated or expensive, and standard SPD assumptions on the near kernel of $A$ can change significantly. Recent analysis developed conditions on $R$ and $P$ to ensure that the nonsymmetric projection is orthogonal in an arbitrary norm \cite{southworth2023compatible}, but relaxation has not yet been incorporated into this theory.

Nonetheless, many important problems are nonsymmetric, motivating us to study this general setting. Consider the interpolation matrix in the following classical AMG form
\begin{equation}\label{eq:classP}
P = \left[\begin{matrix}
   W \\ I 
\end{matrix} \right]
\begin{array}{lr}
\} \: F \\ \} \:  C
\end{array} ,
\end{equation}
where $W$ must be sparse (for efficiency) and represents interpolation from $C$-points to $F$-points. In this setup, the identity matrix below $W$ maintains the values at the coarse points by injecting $C$-points from coarse level to fine level. Different AMG approaches offer different guidelines to construct appropriate interpolation operators of the form \eqref{eq:classP}. Among the most notable and widely accepted guidelines are the so-called  \textit{ideal}~\cite{brannick2010compatible,xu2018ideal} and
\textit{optimal}~\cite{xu2017algebraic,brannick2018optimal} forms of interpolation. In this paper, we focus on the idea of optimal interpolation primarily motivated from \cite{brannick2018optimal}. The aim of \textit{optimal} interpolation is to exactly capture the modes that relaxation leaves behind (i.e., the algebraically smooth modes) in $\operatorname{span}(P)$. The derivation of the optimal interpolation matrix $P_{\sharp}$ and its corresponding classical AMG form allows us to gain a deeper understanding of AMG convergence from a theoretical perspective. In this work, we consider a generalization of the optimal interpolation based AMG framework from \cite{brannick2018optimal} to the nonsymmetric and indefinite setting.

The arrangement of this paper is as follows. \Cref{sec:relaxations} discusses error propagators of different relaxations and a matrix form derivation of the Kaczmarz relaxation method \cite{karczmarz1937angenaherte,popa2008algebraic,brandt1986algebraic} and the corresponding form of symmetrized relaxation. 
\Cref{sec:th_review} reviews the existing two-grid convergence theory \cite{falgout2004generalizing,falgout2005two} and the explicit form of optimal interpolation \cite{brannick2018optimal} which is the main foundation of the proposed theory in this paper. \Cref{sec:orthogonality} introduces a novel linear algebra result regarding a certain matrix induced orthogonality of left and right eigenvectors of a generalized eigenvalue problem (EVP). This result is a key component of the proposed generalized theory for nonsymmetric problems. \Cref{sec:main} provides the two-grid optimal convergence theory for nonsymmetric and indefinite problems. To verify the proposed theory, we present numerical examples in \Cref{sec:num} and the conclusion follows in \Cref{sec:conc}.

\section{Relaxations in matrix forms}\label{sec:relaxations}One of the challenges of convergence theory in the nonsymmetric as well as symmetric indefinite setting is that standard pointwise relaxations such as Jacobi or Gauss-Seidel are not generally guaranteed to be convergent. The classical Kaczmarz method was first introduced in \cite{karczmarz1937angenaherte} as an iterative algorithm for solving a consistent (nonsymmetric or indefinite) system \eqref{eq:main}. Kaczmarz relaxation has been historically used in two-level classical AMG convergence theory \cite{brandt1986algebraic} in part because Kaczmarz is always convergent for arbitrary consistent systems. Convergence of a relaxation method is an essential requirement of the optimal interpolation theory, motivating our primary choice of Kaczmarz. That being said, we will also consider more practical choices of Jacobi and Gauss-Seidel and demonstrate that our theory can be applied there as well. Given our linear system in (\ref{eq:main}), we consider the matrix splitting $A=D+L+U$ where $D$ is the diagonal, $L$ is the strictly lower triangular, and $U$ is the strictly upper triangular part of $A$. In the Jacobi method, the relaxation error propagator is $I -D^{-1}A$. Therefore, for the Jacobi method, we have $M=D$. On the other hand, in the Gauss-Seidel method, the relaxation error propagator is given by $I-(D+L)^{-1}A$.  Therefore, for the Gauss-Seidel method, we have $M=D+L$. 

For Kaczmarz relaxation, we define a new unknown vector $\mathbf{q}$ such that
\begin{equation}\label{transform}
    A^{*}\mathbf{q}=\mathbf{x},
\end{equation}
where $A^{*}$ is the conjugate transpose of $A$. Using \eqref{eq:main}, the equation obtained for $\mathbf{q}$ is 
\begin{equation}\label{kaczmarz}
    AA^{*}\mathbf{q}=\mathbf{b}.
\end{equation}
Kaczmarz is basically Gauss-Seidel iterations on the linear system \eqref{kaczmarz}, where $\mathbf{x} = A^{*} \mathbf{q}$. Note that, we emphasize using $AA^{*}$ instead of $A^{*}A$ based on \cite{brandt1986multigrid}. Let us consider $\widetilde{A}=AA^{*}$. Then the system becomes, $\widetilde{A} \mathbf{q}= \mathbf{b}$. We now decompose $\widetilde{A}$  as $\widetilde{A}=\widetilde{D}+\widetilde{L}+\widetilde{U}$ where $\widetilde{D}$ is the diagonal, $\widetilde{L}$ is the strictly lower triangular, and $\widetilde{U}$ is the strictly upper triangular part of $\widetilde{A}$. Then using Gauss-Seidel,  the solution to  $\widetilde{A} \mathbf{q} = \mathbf{b}$ can be obtained iteratively via $$ \mathbf{q}^{(k+1)} = (\widetilde{D}+\widetilde{L})^{-1} (\mathbf{b} - \widetilde{U} \mathbf{q}^{(k)}). $$
This equation can be further modified as:
$$\mathbf{q}^{(k+1)} =  [I-(\widetilde{D}+\widetilde{L})^{-1}\widetilde{A}]\mathbf{q}^{(k)}+(\widetilde{D}+\widetilde{L})^{-1} \mathbf{b}.$$
Then using $\widetilde{A}=AA^{*}$ and \eqref{transform} we obtain,
\begin{align}\label{kaczmarz_1}
    \mathbf{q}^{(k+1)} &=  \mathbf{q}^{(k)}-(\widetilde{D}+\widetilde{L})^{-1}AA^{*}\mathbf{q}^{(k)}+(\widetilde{D}+\widetilde{L})^{-1} \mathbf{b}\notag\\
       A^{*}\mathbf{q}^{(k+1)} &=  A^{*}\mathbf{q}^{(k)}-A^{*}(\widetilde{D}+\widetilde{L})^{-1}AA^{*}\mathbf{q}^{(k)}+A^{*}(\widetilde{D}+\widetilde{L})^{-1} \mathbf{b}\notag\\
         \mathbf{x}^{(k+1)} &=  [I-(NA^{-*})^{-1}A]\mathbf{x}^{(k)}+(NA^{-*})^{-1} \mathbf{b}\notag\\
        \mathbf{x}^{(k+1)} &=  [I-M^{-1}A]\mathbf{x}^{(k)}+M^{-1} \mathbf{b}
\end{align}
where $N=\widetilde{D}+\widetilde{L}$ and $M=NA^{-*}$. Generally speaking, depending on the system matrix $A$, relaxation $M$ may be convergent or divergent. If $A$ is Hermitian positive definite (HPD), it can be shown that $M$ is a convergent relaxation if and only if $||I-M^{-1}A||_{A}<1$ or $M^{*}+M-A$ is HPD \cite{huang2022learning}. 

The following two operators defined in \cite{falgout2005two} are related to the relaxation $M$:
\begin{align}
    \widetilde{M}&=M^{*}(M^{*}+M-A)^{-1}M,\label{sym_smoother1}\\ 
    \overline{M}&= M (M^{*}+M-A)^{-1}M^{*}.
\end{align}
Note that $I-\widetilde{M}^{-1}A=(I-M^{-1}A)(I-M^{-*}A)$, therefore $\widetilde{M}$ is the symmetrized version of $M$. Similar properties holds for $\overline{M}$. If $M$ is symmetric then $\widetilde{M}=\overline{M}$. In general, $\widetilde{M}$ is used for error analysis and $\overline{M}$ is used in the definition of preconditioners \cite{falgout2005two}. For symmetric $A$, the Galerkin coarse grid operator $A_{c}=P^{*}AP$ ideally inherits the symmetry properties of $A$. This symmetry is crucial for ensuring that the error propagation operators on different levels maintain consistent behavior, which in turn affects the convergence rates. If $A$ is Hermitian but a non-Hermitian relaxation is used, the ``optimal'' transfer operators derived in this paper suggest a nonsymmetric Petrov-Galerkin coarse grid. For this reason, moving forward we will utilize $\widetilde{M}$ rather than $M$ whenever $A$ is Hermitian.

\section{Algebraic two-grid HPD theory review}\label{sec:th_review} In order to conceptualize and establish the generalized optimal convergence theory, we begin with the classical HPD theory \cite{falgout2005two}, which provides the needed concepts of optimal interpolation and of generalized coarse and fine variables. Assuming the system matrix $A$ is HPD, two-grid
convergence for multigrid can be summarized as follows~\footnote{We follow most closely the notation from \cite{brannick2018optimal}, with equation (\ref{kappa_eq}) corresponding to equation (2.1) and the further results (2.8) and (2.9) from~\cite{brannick2018optimal}, and equation (\ref{kappa_sharp_eq}) corresponding to equation (2.1) and the further result (2.3) from~\cite{brannick2018optimal}.}~\cite{falgout2004generalizing,falgout2005two,brannick2018optimal}:
\begin{align}||E_{TG}||_{A}^{2}&\leq 1-\displaystyle\frac{1}{\kappa},\,\,\,\text{where}\,\,\,\kappa\coloneqq\kappa(P\breve{R})=\displaystyle\sup_{\mathbf{e}\ne 0}\frac{||(I-P\breve{R})\mathbf{e}||_{\widetilde{M}}^{2}}{||\mathbf{e}||_{A}^{2}}\geq 1,\label{kappa_eq}\\
    ||E_{TG}||_{A}^{2}&=1-\displaystyle\frac{1}{\kappa_{\sharp}},\,\text{where}\,\kappa_{\sharp}\coloneqq\kappa(\Pi_{\widetilde{M}})=
    \displaystyle\sup_{\mathbf{e}\ne 0}\frac{||(I-\Pi_{\widetilde{M}})\mathbf{e}||_{\widetilde{M}}^{2}}{||\mathbf{e}||_{A}^{2}} \ge 1
 \label{kappa_sharp_eq} \\\nonumber
    &\quad\;\; \mbox{and }\, \Pi_{\widetilde{M}}=P(P^{*}\widetilde{M}P)^{-1}P^{*}\widetilde{M}.
\end{align}

Here,  $\|E_{TG}\|_{A}$ represents the asymptotic convergence factor in the $A$-norm, defined as $\|\mathbf{x}\|_{A}^{2} = (A\mathbf{x}, \mathbf{x})$, where $(\cdot, \cdot)$ is the standard inner product in $\mathbb{C}^{n}$. The symmetrized smoother $\widetilde{M}$ is HPD, and the $\widetilde{M}$-norm is defined as $\|\mathbf{x}\|_{\widetilde{M}}^{2} = (\widetilde{M}\mathbf{x}, \mathbf{x})$. The interpolation matrix $P:\mathbb{C}^{n_{c}}\mapsto\mathbb{C}^{n}$ again maps coarse-grid variables to the fine grid, while ${\breve{R}}:\mathbb{C}^{n}\mapsto\mathbb{C}^{n_{c}}$ (distinct from restriction $R$ in Petrov-Galerkin AMG) is any matrix for which $\breve{R}P=I_{n_{c}}$, the identity on $\mathbb{C}^{n_{c}}$, so that $P\breve{R}$ is a projection onto $\operatorname{range}(P)$. 

The parameter $\kappa$ in \eqref{kappa_eq} quantifies the effectiveness of the coarse-grid correction process. It measures how well the coarse-grid space and interpolation matrix $P$ approximate the fine-grid problem. Smaller values of $\kappa$ (approaching 1) indicate a more effective coarse-grid correction, leading to improved approximation and faster convergence.
The parameter $\kappa_{\sharp}$ in \eqref{kappa_sharp_eq} is a sharp metric that yields the exact two-grid convergence rate and is based on the projection $\Pi_{\widetilde{M}}$, i.e., the $\widetilde{M}$-orthogonal projection onto $\operatorname{range}(P)$. 
If we minimize $\kappa_{\sharp}$ in (\ref{kappa_sharp_eq}) (see equation \eqref{eq:etg_bigK_detail} in Lemma \ref{lem:optimal}), then we find the optimal, or best $P$ for two-grid convergence for a given $n_c$.
In other words, this yields \textit{optimal} AMG interpolation \cite{brannick2018optimal,xu2017algebraic}, which is our primary interest. 
If we instead minimize with respect to \eqref{kappa_eq} and replace the $\widetilde{M}$ norm with $\ell_{2}$-norm when defining $\kappa$, this yields the standard \textit{ideal} AMG
interpolation \cite{falgout2004generalizing}.\footnote{$
P_{\operatorname{ideal}}=\begin{bmatrix}
 -A_{f\!f}^{-1}A_{f\!c}\\
 I
\end{bmatrix}
$, which among other perspectives, constructs a Schur complement coarse grid. For details on this derivation, see equations (2.8), (2.9) and Lemma 2 in~\cite{brannick2018optimal}.} 
Some key assumptions in this theory:
\begin{itemize}
    \item $A \in \mathbb{C}^{n \times n}$ is an  HPD matrix.
    \item Convergence is considered in the well defined energy-norm or $A$-norm.
    \item Reliance on  $\Pi_{\widetilde{M}}=P(P^{*}\widetilde{M}P)^{-1}P^{*}\widetilde{M}$, i.e., the $\left( \cdot, \cdot \right)_{\widetilde{M}}$ orthogonal projection onto $\operatorname{range}(P)$.
\end{itemize}

We now outline the recently derived explicit form of optimal interpolation \cite{brannick2018optimal} for the HPD case, namely the interpolation operator that minimizes the convergence rate of the two-grid method over all possible choices of $P$, with respect to an HPD relaxation operator $\widetilde{M}$ and coarse-grid dimensionality $n_c$.
\begin{lemma}\label{lem:optimal}
\noindent Given $n_{c}$ and $\widetilde{M}$, let $\mathbf{y}_{i}$ solve the generalized EVP 
\begin{equation} \label{gen_prob_first}
A\mathbf{y}_{i}=\mu_{i}\widetilde{M}\mathbf{y}_{i}  
\end{equation}
where $i=1,2,\dots,n$; 
$P:\mathbb{C}^{n_{c}}\to\mathbb{C}^{n}$ is full rank, $\mu_{1}\leq \mu_{2} \leq \cdots \leq\mu_{n_{c}}\leq \cdots  \leq \mu_{n}$ denote the eigenvalues and $\mathbf{y}_{i}$ denote the corresponding $\widetilde{M}$-orthonormal eigenvectors, i.e., $\left(\mathbf{y}_{i}, \mathbf{y}_{j}\right)_{\widetilde{M}}=\mathbf{y}_{i}^{*}\widetilde{M}\mathbf{y}_{j} = \delta_{ij}$.
Then the minimal convergence rate of the two-grid method over all $P$ is
attained by the optimal interpolation $P_{\sharp}$, where
\begin{equation}\label{eq:etg_bigK}
    ||E_{TG}||_{A}^{2}=  1 - \displaystyle\frac{1}{\mathcal{K}},  
\end{equation}
\begin{equation} \label{eq:etg_bigK_detail}
    \mathcal{K} :=  \displaystyle{\inf_{\operatorname{dim}(\operatorname{range}(P))\, =\, n_{c}} \kappa_{\sharp}(P)}
    \,=\, \displaystyle{\inf_{\operatorname{dim}(\operatorname{range}(P))\, =\, n_{c}}} \,\, \sup_{\mathbf{e}\ne 0}\frac{||(I-\Pi_{\widetilde{M}})\mathbf{e}||_{\widetilde{M}}^{2}}{||\mathbf{e}||_{A}^{2}}
    = \displaystyle \frac{1}{\mu_{n_{c}+1}},
\end{equation}
implying $ \mathcal{K}=\displaystyle \frac{1}{\mu_{n_{c}+1}}$. Thus, substituting for $\mathcal{K}$ in equation \eqref{eq:etg_bigK} yields
\begin{equation}\label{eq:etg_A_norm}
    ||E_{TG}||_{A}^{2}=1 -\mu_{n_{c}+1},
\end{equation}
 where the optimal interpolation operator $P_{\sharp}$ 
satisfies
     \begin{equation}\label{con_rate}
         \operatorname{range}(P_{\sharp}) = \operatorname{range}\left( 
\begin{pmatrix} \mathbf{y}_{1} & \mathbf{y}_{2} & \cdots & \mathbf{y}_{n_{c}} \end{pmatrix}\right).
     \end{equation}
\end{lemma}
The optimal convergence rate for the two-grid method can thus be found by selecting any interpolation with a range equivalent to $P_{\sharp}$, that is, by setting $P=P_{\sharp} Z$, for coarse-grid change of basis matrix $Z$, where we assume the existence of $Z^{-1}$. The detailed proof of the above lemma can be found in \cite{xu2017algebraic,brannick2018optimal} where it uses the definition of the energy norm and Courant-Fischer Min-max representation.

The fundamental concept underlying \Cref{lem:optimal} is that the convergence rate of the two-grid method is minimized by selecting an optimal interpolation operator $P_{\sharp}$ that maps onto the subspace spanned by the first $n_c$ eigenvectors of the generalized eigenvalue problem \eqref{gen_prob_first}. The convergence rate is determined by the $A$-norm of the two-grid error transfer operator, which is expressed as \eqref{eq:etg_bigK}, where $\mathcal{K}$ is optimized over all possible interpolation operators. The critical eigenvalue $\mu_{n_c+1}$ governs this rate, as the expression (\ref{eq:etg_bigK_detail}) simplifies to $\mathcal{K} = \frac{1}{\mu_{n_{c}+1}}$. Consequently, the convergence rate is minimized with value in \eqref{eq:etg_A_norm}. For the proof, see~\cite{brannick2018optimal}.  One key implication is that a larger $\mu_{n_{c}+1}$ corresponds to a faster optimal convergence rate.

\section{Induced orthogonality of generalized eigenvectors}\label{sec:orthogonality}

The classical AMG theory for SPD matrices is built around orthogonal decompositions of error. One of the major difficulties of nonsymmetric theory is in forming such an orthogonal decomposition, particularly separating coarse-grid correction and relaxation. Here we establish a linear algebra result for general non-Hermitian matrices with diagonalizable product, where left and right generalized eigenvectors enjoy a matrix-induced orthogonality with respect to each matrix of the matrix pencil.\footnote{Noting that the set of nondiagonalizable matrices in $\mathbb{C}^{n\times n}$ is measure zero, we point out diagonalizability is a fairly weak assumption.} Although for non-Hermitian matrices, there is no corresponding sense of inner product, the induced orthogonality still proves critical in our later analysis of error propagation. 

\begin{lemma}\label{lem:orth}
    Let $A,M$ be square matrices such that $M$ is invertible and $M^{-1}A$ is diagonalizable. Consider the left and right generalized eigenvectors, $V_l$ and $V_r$, respectively, of the matrix pencil $(A,M)$, defined such that
    \begin{subequations}\label{eq:gep}
    \begin{align}
        AV_r & = MV_r\Lambda, \\
       \ V_l^*A & = \Lambda V_l^*M,\label{eq:gep-left}
    \end{align}
    \end{subequations}
    where $\Lambda$ is a diagonal matrix of eigenvalues. Then $V_l$ and $V_r$ induce a matrix-based orthogonality, satisfying
    \begin{subequations}
    \begin{align}
       {V_l}^*AV_r &= D_a, \\
       {V_l}^*MV_r &= D_m,
    \end{align}
    \end{subequations}
    for diagonal matrices $D_a,D_m$.
  
\end{lemma}
\begin{proof}

To prove the matrix-based orthogonality relations, we first transform the generalized EVPs into equivalent standard eigenvalue problems using the invertibility of $M$. Define the matrices $M^{-1}A$ and $AM^{-1}$, which lead to the following eigenvalue problems:
\begin{subequations}
\begin{align}
    M^{-1}A \hat{V}_{1,r} &= \hat{V}_{1,r} \Lambda, \quad &\text{(right eigenvectors of \(M^{-1}A\))}, \label{eq:right1}\\
    A^*M^{-*} \hat{V}_{1,l} &= \hat{V}_{1,l} \overline{\Lambda}, \quad &\text{(left eigenvectors of \(M^{-1}A\))}, \label{eq:left1}\\
    AM^{-1} \hat{V}_{2,r} &= \hat{V}_{2,r} \Lambda, \quad &\text{(right eigenvectors of \(AM^{-1}\))}, \label{eq:right2}\\
    M^{-*}A^* \hat{V}_{2,l} &= \hat{V}_{2,l} \overline{\Lambda}, \quad &\text{(left eigenvectors of \(AM^{-1}\))}. \label{eq:left2}
\end{align}
\end{subequations}
Here, we use $\hat{V}_{1,r}$, $\hat{V}_{1,l}$, $\hat{V}_{2,r}$, and $\hat{V}_{2,l}$ to denote standard eigenvectors, with subscripts $l$ and $r$ for left and right eigenvectors, and subscripts $1,2$ for the operators $M^{-1}A$ and $AM^{-1}$, respectively.  The equivalence of eigenproblems between $M^{-1}A$ and $AM^{-1}$ allows us to relate eigenvectors. Also, the notation $\overline{\Lambda}$ denotes the element-wise complex conjugate of $\Lambda$. Taking conjugate transpose $(\cdot)^{*}$ of \eqref{eq:left2} and multiplying on the left and right by $\hat{V}_{2,l}^{-*}$, we find
\begin{equation*}
     AM^{-1}\hat{V}_{2,l}^{-*} =\hat{V}_{2,l}^{-*}\Lambda.
\end{equation*} Comparing this result  with \eqref{eq:right2} yields the relationship
\begin{equation}\label{eq:rel1}
    \hat{V}_{2,r} = \hat{V}_{2,l}^{-*} D_{1},
\end{equation}
where the diagonal matrix $D_1$ accounts for the scaling required to match the normalization of left and right eigenvectors.
Next, multiplying \eqref{eq:right1} by $A$ on the left, 
\begin{equation}\label{eq:compare_D_2}
    (AM^{-1})(A \hat{V}_{1,r}) = (A\hat{V}_{1,r} )\Lambda.
\end{equation}
Comparing this equation with \eqref{eq:right2}, we find \begin{equation}\label{eq:rel2}
    \hat{V}_{2,r} = A \hat{V}_{1,r} D_2,
\end{equation}
where $D_2$ is a diagonal scaling matrix which adjusts the contribution of $A$ applied to the right eigenvectors of $M^{-1}A$.
Similarly, multiplying \eqref{eq:right1} by $M$ on the left and inserting the identity, we find
 \begin{equation*}
     A(M^{-1}M)\hat{V}_{1,r}=M\hat{V}_{1,r}\Lambda,
 \end{equation*}
which further yields
\begin{equation}\label{eq:M_invertibility}
    (AM^{-1})(M\hat{V}_{1,r})=(M\hat{V}_{1,r}) \Lambda.
\end{equation}
Comparing this equation with \eqref{eq:right2} leads to
\begin{equation}\label{eq:rel3}
    \hat{V}_{2,r} = M \hat{V}_{1,r} D_3, 
\end{equation}
where $D_{3}$ is a diagonal scaling matrix which adjusts the normalization when scaling by $M$ instead of $A$. 
Combining \eqref{eq:rel1} with \eqref{eq:rel2} and \eqref{eq:rel1} with \eqref{eq:rel3} yield
\begin{equation}
    \hat{V}_{2,l}^*A\hat{V}_{1,r} = D_a, \quad \hat{V}_{2,l}^*M\hat{V}_{1,r} = D_m, \label{eq:orthogonality}
\end{equation}
where $D_a$ and $D_m$ are diagonal matrices representing the orthogonality structure induced by $A$ and $M$, respectively. Furthermore, $D_1$, $D_2$, and $D_3$ are related to $D_a$ and $D_m$ as follows:
\begin{equation*}
    D_a = D_{1}D_2^{-1}, \hspace{4ex}
    D_m = D_1 D_3^{-1}.
\end{equation*}
Finally, note that $\hat{V}_{1,r}$ and $\hat{V}_{2,l}$ correspond to the generalized eigenvectors $V_r$ and ${V_l}$ of the original generalized EVP. Specifically, $V_r = \hat{V}_{1,r}$ and $V_l = \hat{V}_{2,l}$. Using this identification, we rewrite the orthogonality relations 
\begin{equation}\label{eq:final_orthogonality}
  {V_l}^*AV_r = D_{a}, \quad {V_l}^*MV_r = D_{m}, 
\end{equation}which completes the proof.
\end{proof}

For Hermitian indefinite $A$ and $M$, we can then define an indefinite bilinear form, e.g. $[\cdot,\cdot]_{M}=(M\mathbf{x},\mathbf{y})=\mathbf{y}^{*}M\mathbf{x}$ \cite{gohberg2006indefinite,kahl2009adaptive}, where $(\cdot, \cdot)$ denotes the standard inner product in $\mathbb{C}^{n}$ and  $\mathbf{x},\mathbf{y}\in\mathbb{C}^{n}$. Due to the fact that $M$ is indefinite, this is not an inner product in the usual sense, but it provides further structure to the eigenvector orthogonality. Noting that for real-valued Hermitian $A$ and $M$ we have $V_l = \overline{V}_r$ in \Cref{lem:orth} (take the adjoint and conjugate of \eqref{eq:gep-left}), we arrive at the following Corollary.

\begin{corollary}\label{cor:hermitian_orthogonality}
    Let $A$ and $M$ be real-valued Hermitian matrices and consider the generalized eigenvectors $V$ satisfying
    \begin{equation}
         AV = MV\Lambda,
    \end{equation}
    for eigenvalues $\Lambda$. Then,
    \begin{equation}
        V^TAV = \overline{V}^*AV = D_a, \hspace{4ex} V^TMV = \overline{V}^*MV = D_m,
    \end{equation}
    for diagonal matrices $D_a,D_m$.\footnote{Note that $\bar{(\cdot)}$ represents complex conjugation, while $(\cdot)^*$ denotes the conjugate transpose. Consequently, the operation $\bar{(\cdot)}^*$ simplifies to the strict (non-conjugate) transpose, which is denoted by $(\cdot)^T$.  Therefore, $\overline{V}^* = V^T$.}
\end{corollary}

Equivalently to the above corollary, in indefinite inner products over $A$ and $M$, generalized eigenvectors $V$ form a permuted orthogonal basis:
    \begin{equation*}
        [\mathbf{v}_{\lambda_i},\mathbf{v}_{\lambda_j}]_{A} = \begin{cases} 
        (D_a)_{ii} & \text{if} \,\,\,\lambda_i=\overline{\lambda_j} \\
        0 & \text{otherwise} 
        \end{cases}, 
        \hspace{4ex}
        [\mathbf{v}_{\lambda_i},\mathbf{v}_{\lambda_j}]_{M} = \begin{cases} 
        (D_m)_{ii} & \text{if} \,\,\,\lambda_i=\overline{\lambda_j} \\
        0 & \text{otherwise}
        \end{cases},
    \end{equation*}
where complex conjugate eigenvector pairs form an orthogonal subspace of dimension two, where they are orthogonal to themselves in the $A$-and $M$-inner products.

\section{Optimal AMG theory for non-HPD problems}\label{sec:main}
In the context of nonsymmetric AMG, coarse grid correction typically involves a nonorthogonal projection in any known inner product. Construction of coarse-grid correction in the nonsymmetric setting generally relies on constructing a separate restriction and interpolation operator $R$ and $P$ respectively in a compatible sense rather than using $R=P^{*}$ (see \cite{southworth2023compatible}). Consider the two grid error transfer operator with only a pre-relaxation step with some $M\approx A$ given by \begin{equation}\label{ETG_nonsym}
    E_{TG}(P,R)=(I-P(RAP)^{-1}RA)(I-M^{-1}A).
\end{equation}
We now propose the two-grid optimal convergence theory for nonsymmetric problems as follows:
\begin{theorem}\label{th:optimalPR_nonsymmetric}
Let $A\in\mathbb{C}^{n\times n}$ be a general non-singular matrix. Given a subset   $\mathcal{I}\subseteq \{1,2,\dots,n\}$ with cardinality $|\mathcal{I}|=n_{c}$, and non-singular smoother $M$ such that $M^{-1}A$ is diagonalizable, Consider the left and right generalized eigenvectors, $V_l=\left[\mathbf{\check{v}_{i}}\right]_{i=1}^{n}$ and $V_r=\left[\mathbf{\hat{v}_{i}}\right]_{i=1}^{n}$, of the matrix pencil $(A,M)$, respectively defined as in \eqref{eq:gep} where the corresponding eigenvalues are  $\left\{\lambda_{i}\right\}_{i=1}^{n}$.  Define the interpolation operator $P$ and restriction operator $R$ such that their columns are subsets of the right and left generalized eigenvectors, respectively:
\begin{align*}
\operatorname{range}(P) &= \operatorname{range}\left(
\begin{pmatrix} \mathbf{\hat{v}}_{\ell}:\ell \in \mathcal{I}  \end{pmatrix}\right),\\
\operatorname{range}(R^{*}) &= \operatorname{range}\left(
\begin{pmatrix} \mathbf{\check{v}}_{\ell} :\ell \in \mathcal{I}   \end{pmatrix}\right).
\end{align*}
Then, the spectral radius of the two-grid error transfer operator \eqref{ETG_nonsym} is given by
\begin{equation}\label{identity_nonsymmetric}
    \rho\left(E_{TG}(P,R)\right)=\max_{\ell\notin \mathcal{I}}|1-\lambda_{\ell}| .
\end{equation}
\end{theorem}

\begin{proof}
Here we prove the result for interpolation and restriction operators $P$ and $R$, whose columns are subsets of the right and left generalized eigenvectors. The proof extends directly from the structure of the generalized EVP and the properties of the error transfer operator. 

Consider the generalized EVP for the general non-singular matrix $A$ in the following matrix equation form \begin{equation}\label{main_th_nonsym}
    A\, V_{r}=MV_{r}\Lambda,
\end{equation} where $M$ is the relaxation operator, the columns of $V_{r}=[\mathbf{\hat{v}}_{1},\mathbf{\hat{v}}_{2},...,\mathbf{\hat{v}}_{n}]$ are the right generalized eigenvectors, and the diagonal elements of $\Lambda:=\operatorname{diag}(\lambda_{1},\lambda_{2},...,\lambda_{n})$ are the generalized eigenvalues. Similarly, the left generalized eigenvectors  $V_{l} = [\mathbf{\check{v}}_{1}, \mathbf{\check{v}}_{2}, \ldots, \mathbf{\check{v}}_{n}]$  satisfy:
\begin{equation}\label{main_th_nonsym_left}
    V_{l}^{*} A = {\Lambda}V_{l}^{*} M .
\end{equation}
The error transfer operator of the two-level
method is
\begin{equation}\label{eq: two_grid_error_transfer}
     E_{TG}(P,R)=(I-PA_{c}^{-1}RA)(I-M^{-1}A),
\end{equation} where $A_{c} = RAP$ is the Petrov-Galerkin variational coarse-level matrix.

Define the interpolation and restriction operators as:
\begin{equation}
    P = \begin{bmatrix} \mathbf{\hat{v}}_{\ell} \colon \ell \in \mathcal{I} \end{bmatrix}, \quad
R^{*} = \begin{bmatrix} \mathbf{\check{v}}_{\ell} \colon \ell \in \mathcal{I} \end{bmatrix},
\end{equation}
where  $\mathcal{I} \subseteq \{1, 2, \ldots, n\}$  has cardinality  $|\mathcal{I}| = n_{c}$. Let  $\mathcal{I}^{c} = \{1, 2, \ldots, n\} \setminus \mathcal{I}$  denote the complement of $\mathcal{I}$.
Using \eqref{main_th_nonsym}, the relaxation applied to $V_{r}$ is:
\begin{equation}\label{eq: relax_applied_Vr}
 (I-M^{-1}A)V_{r}= V_{r}-M^{-1}AV_{r}=V_{r}-V_{r}\Lambda= V_{r}(I-\Lambda).  
\end{equation}

Expanding  $V_{r}$  into columns corresponding to  $\mathcal{I}$  and  $\mathcal{I}^{c}$, we write:
\begin{equation}\label{eq: V_r_expanded}
    V_{r} = \begin{bmatrix} P & S \end{bmatrix},
\end{equation}
where $S = \begin{bmatrix} \mathbf{\hat{v}}_{\ell} \colon \ell \in \mathcal{I}^{c} \end{bmatrix}$. Substituting this into \eqref{eq: relax_applied_Vr}, we have:
\begin{equation}\label{eq: substitute_after_split}
(I - M^{-1} A) V_{r} = V_{r} (I - \Lambda) = \begin{bmatrix} P & S \end{bmatrix}
\begin{bmatrix} I_{n_{c}} - \Lambda_{\mathcal{I}} & 0 \\ 0 & I_{n - n_{c}} - \Lambda_{\mathcal{I}^{c}} \end{bmatrix},
\end{equation}
where $\Lambda_{\mathcal{I}} = \operatorname{diag}(\lambda_{\ell} \colon \ell \in \mathcal{I})$  and  $\Lambda_{\mathcal{I}^{c}} = \operatorname{diag}(\lambda_{\ell} \colon \ell \in \mathcal{I}^{c})$.
Since $R$ and $S$ consist of left and right generalized eigenvectors of the matrix pencil $(A,M)$, using \Cref{lem:orth} it follows that \begin{equation}\label{eq:induced_orthogonality}
    RMS=0.
\end{equation}
Again, from \eqref{main_th_nonsym} and \eqref{eq: V_r_expanded} , we have
\begin{equation}
    A \begin{bmatrix} P & S \end{bmatrix} = M \begin{bmatrix} P & S \end{bmatrix} \begin{bmatrix}  \Lambda_{\mathcal{I}} & 0 \\ 0 &  \Lambda_{\mathcal{I}^{c}} \end{bmatrix},
\end{equation}
which implies 
\begin{equation}\label{eq: gep_splitted}
    AP=MP \Lambda_{\mathcal{I}}, \quad AS=MS \Lambda_{\mathcal{I}^{c}}.
\end{equation}
Then, from \eqref{eq: two_grid_error_transfer}, we have
\begin{align*}
 E_{TG}(P,R)V_{r}&=(I-P(RAP)^{-1}RA)(I-M^{-1}A)V_{r}\\
 &=(I-P(RAP)^{-1}RA)V_{r}(I-\Lambda)[\text{using}\eqref{eq: substitute_after_split}]\\
 &=(V_{r}-P(RAP)^{-1}RAV_{r})(I-\Lambda)\\
 &=(V_{r}-P(RAP)^{-1}RMV_{r}\Lambda)(I-\Lambda)[\text{using}\eqref{main_th_nonsym}]\\
 &=(V_{r}-P(RMP \Lambda_{\mathcal{I}})^{-1}RMV_{r}\Lambda)(I-\Lambda)[\text{using}\eqref{eq: gep_splitted}]\\
 &=(V_{r}-P\Lambda_{\mathcal{I}}^{-1}(RMP )^{-1}RMV_{r}\Lambda)(I-\Lambda)\\
 &=\left(\begin{bmatrix} P & S \end{bmatrix}-P\Lambda_{\mathcal{I}}^{-1}(RMP )^{-1}RM\begin{bmatrix} P & S \end{bmatrix}\Lambda\right)(I-\Lambda)[\text{using}\eqref{eq: V_r_expanded}]\\
 &=\left(\begin{bmatrix} P & S \end{bmatrix}-P\Lambda_{\mathcal{I}}^{-1}(RMP )^{-1}\begin{bmatrix} RMP & RMS \end{bmatrix}\Lambda\right)(I-\Lambda)\\
 &=\left(\begin{bmatrix} P & S \end{bmatrix}-P\Lambda_{\mathcal{I}}^{-1}(RMP )^{-1}\begin{bmatrix} RMP & 0 \end{bmatrix}\Lambda\right)(I-\Lambda)[\text{using}\eqref{eq:induced_orthogonality}]\\
 &=\left(\begin{bmatrix} P & S \end{bmatrix}-P\Lambda_{\mathcal{I}}^{-1}\begin{bmatrix} \Lambda_{\mathcal{I}} & 0 \end{bmatrix}\right)(I-\Lambda)\\
 &=\left(\begin{bmatrix} P & S \end{bmatrix}-\begin{bmatrix} P & 0 \end{bmatrix}\right)(I-\Lambda)\\
 &=\begin{bmatrix} 0 & S \end{bmatrix}\begin{bmatrix} I_{n_{c}} - \Lambda_{\mathcal{I}} & 0 \\ 0 & I_{n - n_{c}} - \Lambda_{\mathcal{I}^{c}} \end{bmatrix}\\
 &=V_{r}\begin{bmatrix} 0 & 0 \\ 0 & I_{n - n_{c}} - \Lambda_{\mathcal{I}^{c}} \end{bmatrix}.
 \end{align*}
 Thus,  $E_{TG}(P, R)$  annihilates the components in the range of  $P$, while its eigenvalues on the subspace spanned by  $S$  are  $|1 - \lambda_{\ell}|$  for  $\ell \in \mathcal{I}^{c}$ . The spectral radius is therefore
\begin{equation}\label{eq:identity_nonsymmetric_max}
    \rho(E_{TG}(P, R)) = \max_{\ell \in \mathcal{I}^{c}} |1 - \lambda_{\ell}|=\max_{\ell \notin \mathcal{I}} |1 - \lambda_{\ell}|.
\end{equation}
\end{proof}
One of the most important aspects of this proof is that, without even introducing an $A$-norm, we established the identity \eqref{eq:identity_nonsymmetric_max} using a direct proof approach. Using interpolation and restriction operators whose columns are subset of the right and left generalized eigenvectors, respectively, we were able to demonstrate a two grid convergence theory, even though an energy norm does not exist for a nonsymmetric $A$. It should be noted that on $\mathbb{C}^{n\times n}$, the spectral radius does not define a norm. For $E_{TG}$, the identity \eqref{eq:identity_nonsymmetric_max} is therefore not a $2$-norm bound.

An immediate consequence of the above theorem is the following corollary of pseudo-optimal asymptotic convergence under optimal interpolation $P_\sharp$ and optimal restriction $R_\sharp$, by showing that $P_\sharp$ and $R_\sharp$ minimize the spectral radius of two-grid error propagation over the space of transfer operators with ranges defined by some span of generalized eigenvectors. 
\begin{corollary}[Optimality of transfer operators]\label{cor:optimality_transfer}
    Let $A,M\in\mathbb{C}^{n\times n}$ be such that $M$ is invertible and $M^{-1}A$ is diagonalizable. Consider the inter-grid transfer operators $P, R^* \in \mathbb{C}^{n\times n_c}$ to be defined to span some set of $n_{c}$ right and left generalized eigenvectors of the matrix pencil $(A,M)$, respectively. Assume that the corresponding eigenvalues $\left\{\lambda_{i}\right\}_{i=1}^{n}$ are ordered such that $|1-\lambda_{1}|\geq |1-\lambda_{2}| \geq \cdots  \geq|1-\lambda_{n}|$. Define the optimal interpolation $P_{\sharp}$ and restriction $R_{\sharp}$ to satisfy
    \begin{align*}
        \operatorname{range}(P_{\sharp}) &= \operatorname{range}\left(
        \begin{pmatrix} \mathbf{\hat{v}}_{1} & \mathbf{\hat{v}}_{2} & \cdots & \mathbf{\hat{v}}_{n_{c}} \end{pmatrix}\right),\\
        \operatorname{range}(R_{\sharp}^{*}) &= \operatorname{range}\left(
        \begin{pmatrix} \mathbf{\check{v}}_{1} & \mathbf{\check{v}}_{2} & \cdots & \mathbf{\check{v}}_{n_{c}} \end{pmatrix}\right),
    \end{align*}
    where $\mathbf{\hat{v}}_{i}$ and $\mathbf{\check{v}}_{i}$ are the right and left eigenvectors corresponding to the eigenvalue $\lambda_{i}$, respectively. Then, over this space of transfer operators, the optimal interpolation $P_\sharp$ and optimal restriction $R_\sharp$ minimize the spectral radius of the two-grid method, i.e.,
    \begin{equation}\label{eq: cor_iden}
        \min_{P, R^{*} \in \mathbb{C}^{n\times n_c}} \rho(E_{TG}(P, R)) = 
    \rho(E_{TG}(P_\sharp,R_\sharp)) = 
    |1-\lambda_{n_{c}+1}|.
    \end{equation}
\end{corollary}
\begin{proof}
By assumption, the transfer operators $P$ and $R^{*}$ correspond to the span of some set of $n_c$ right and left generalized eigenvectors of the matrix pencil $(A,M)$, respectively. Due to the fact that $A$ and $M$ are general non-singular matrices, complex conjugate eigenvalues can occur in the spectrum. Hence, the eigenvalues of the generalized EVP are analyzed in the context of the relaxation operator $I-M^{-1}A$. To facilitate this, we sort the eigenvalues based on the function \( |1 - \lambda|^2 \). For a given eigenvalue $\lambda = re^{i\theta}$, this function can be expressed as  $|1 - re^{i\theta}|^2$, which simplifies to  $1 + r^2 - 2r\cos{\theta}$. Accordingly, we arrange the eigenvalues in descending order based on the values of  $1 + r^2 - 2r\cos{\theta}$, effectively sorting them by their squared distance from $1$ in the complex plane. Any eigenvalue for which $|1-\lambda|>1$ corresponds to a divergent mode. Therefore, the spectral radius of the error transfer operator $ E_{TG}(P, R) $ is minimized when the maximum value of $ |1 - \lambda| $ is minimized. If we choose arbitrary $P$ and $R^{*}$ that span a different set of $n_{c}$ generalized eigenvectors than $P_{\sharp}$ and $R_{\sharp}$, the excluded eigenvalues will have a larger deviation from  $1$ (i.e., larger $|1 - \lambda|$) than the excluded eigenvalue with largest deviation from 1 in the optimal case, i.e., $\lambda_{n_{c}+1}$. Thus only $P_\sharp$ and $R_\sharp$ can achieve the minimum possible spectral radius for the error transfer operator, proving our corollary.

\end{proof}

We can specialize the optimal two-grid transfer operators from \Cref{cor:optimality_transfer} to real-valued Hermitian (potentially indefinite) operators, such as those arising from saddle point systems, which will allow a natural Galerkin $R = P^*$ set of optimal transfer operators.

\begin{corollary}[Real-valued Hermitian operators]\label{cor:hermitian_indefnite}
Let $A\in\mathbb{R}^{n\times n}$ be a Hermitian matrix. Given $n_{c}$ and a non-singular Hermitian relaxation operator $M\in\mathbb{R}^{n\times n}$ such that $M^{-1}A$ is diagonalizable, consider the generalized EVP
 \begin{equation}\label{gen_prob_hermitian}
A\mathbf{v}_{i}=\lambda_{i}M\mathbf{v}_{i}, 
 \end{equation}
where $i=1,2,\dots,n$; and the corresponding eigenvalues $\left\{\lambda_{i}\right\}_{i=1}^{n}$ are ordered such that $|1-\lambda_{1}|\geq |1-\lambda_{2}| \geq \cdots  \geq|1-\lambda_{n}|$. Define the optimal interpolation operator $P_{\sharp}$ to satisfy
\begin{equation*}
\operatorname{range}(P_{\sharp}) = \operatorname{range}\left(
\begin{pmatrix} \mathbf{v}_{1} & \mathbf{v}_{2} & \cdots & \mathbf{v}_{n_{c}} \end{pmatrix}\right),
\end{equation*}
where we assume that when complex conjugate eigenpairs are considered, if one complex conjugate eigenvector is placed in $\operatorname{range}(P_{\sharp})$, then the other must be at well.\footnote{This restriction guarantees invertibility of  $P_{\sharp}^{*}MP_{\sharp}$; see Remark \ref{complex_remark} for more details.} Then letting $R_{\sharp}=P_{\sharp}^{*}$, the spectral radius of the two-grid error transfer operator is given by
\begin{equation}\label{identity}
\rho\left(E_{TG}(P_{\sharp})\right)= |1 - \lambda_{n_{c}+1}|.
\end{equation}
\end{corollary}
\begin{proof}
    Let $[\cdot,\cdot]_{M} = (M\,\,\cdot,\cdot)$ be an indefinite inner product \cite{gohberg2006indefinite} with corresponding
Hermitian  matrix $M$. Let the columns of $P_{\sharp}$ be given by $
\begin{pmatrix}\mathbf{v}_{1} & \mathbf{v}_{2}  &\cdots &\mathbf{v}_{n_c}\end{pmatrix}$ and define the $M$-orthogonal complement (w.r.t $[\cdot,\cdot]_{M}$) as $\mathcal{S}_{\sharp}$ with columns given by $ 
\begin{pmatrix} \mathbf{v}_{n_{c}+1} & \mathbf{v}_{n_{c}+2} & \cdots & \mathbf{v}_{n} \end{pmatrix}.$ We use specific eigenvectors as columns for convenience, but note that because coarse-grid correction is invariant to a change of basis in $P$, and due to the $M$-orthogonal decomposition of eigenvectors, the results naturally generalize to $P_\sharp$ and $S_\sharp$ with ranges given by the specified eigenvectors. As a direct consequence of \Cref{cor:hermitian_orthogonality}, $P_\sharp^{*}MS_\sharp=0$. Therefore, using $R_{\sharp}=P_{\sharp}^{*}$ in \Cref{cor:optimality_transfer}, it immediately follows that $\rho\left(E_{TG}(P_{\sharp})\right)= |1 - \lambda_{n_{c}+1}|$.
\end{proof}

\begin{corollary}[Necessary conditions]
    Let the transfer operators $P$ and $R^*$ be defined to span some set of $n_c$ corresponding right and left generalized eigenvectors, respectively, with indices $\mathcal{I} = \{i_1,..,i_{n_c}\}$. Then, in any consistent matrix norm, $\|\cdot\|$, a necessary condition for convergence after $p>1$ iterations, $\|E_{TG}^p\| < 1$, is that
    \begin{equation*}
        | 1 - \lambda_\ell | < 1
    \end{equation*}
    for all $\ell\in\mathcal{I}$, where $\lambda_\ell$ is the $\ell$th generalized eigenvalue.
\end{corollary}
\begin{proof}
    Recall that for matrix $A$, integer $k$, and consistent matrix norm, $\|\cdot\|$, we have $\rho(A) \leq \|A^k\|^{1/k}$. In addition, $\|A^k\| < 1$ if and only if $\|A^k\|^{1/k} < 1$. Thus $\rho(E_{TG})<1$ is a necessary condition for convergence in norm, $\|E_{TG} \| < 1$. From \Cref{th:optimalPR_nonsymmetric} it follows immediately that $|1-\lambda_{l}|<1$ for all $l\in\mathcal{I}$ is a necessary condition for $\|E_{TG} \| < 1$. 
\end{proof}
\begin{remark}\label{divergence_remark}
The identity \eqref{eq: cor_iden} indicates the two-grid method's asymptotic convergence (or divergence) rate. In the case of convergent relaxation, e.g., Kaczmarz, where 
\begin{equation*}\rho (I-M^{-1}A)=\max_{i}|1-\lambda_{i}(M^{-1}A)|<1,
\end{equation*}
this implies $0< \lambda_{i}(M^{-1}A)<2$.  That is, all eigenvalues are less than 2 in magnitude and the quantity $|1-\lambda_{n_{c}+1}|<1$, indicating convergence of the overall method. However, we establish the theory in a way that also applies to the case when the relaxation $I-M^{-1}A$ is divergent, i.e., if $\rho(I-M^{-1}A)>1$.  In this case, the two-grid method overall could diverge or converge, depending on the coarse-grid correction.  If the coarse grid correction suppresses the ``bad" eigenmodes to a sufficient degree, then we can still obtain convergence in this case. More specifically, if the coarse-grid correction takes $n_c$ large enough to remove all divergent eigenmodes of the relaxation, then the overall method will still be convergent. One particular case of divergence based on block Jacobi supporting this statement is briefly discussed in \cref{sec:adv-diff}. 
\end{remark}
\begin{remark}\label{complex_remark}
     The generalized EVP  $A\, \mathbf{x}=\lambda M\mathbf{x},$ in \Cref{cor:hermitian_indefnite}
 consists of two real-valued Hermitian matrices $A$ and $M$. In the classical AMG theory for HPD problems, $M$ is positive definite, hence the problem can be reduced to a Hermitian eigenvalue of the form $M^{-1/2}AM^{-1/2}\mathbf{y}=\lambda \mathbf{y}$ where $M^{1/2}$ is the positive definite square root of $M$. Therefore the generalized eigenvalues are real \cite{stewart1979pertubation}. But in the proposed generalization of the optimal convergence theory for Hermitian indefinite operators, we may encounter complex conjugate eigenvalues of the generalized EVP \eqref{gen_prob_hermitian}. Inclusion of eigenvectors in the span of optimal interpolation $P_{\sharp}$ corresponding to \emph{both} complex conjugate eigenvalues is necessary to ensure that $P_{\sharp}^{*}MP_{\sharp}$ is invertible and $P_{\sharp}^{*}MS_{\sharp} =0$. This corresponds to the conjugate pair orthogonality discussed in \Cref{cor:hermitian_orthogonality}, and can be seen by considering an illustrative example where we take two pairs of complex conjugate eigenvalues (in total four eigenvalues) and their corresponding eigenvectors as the columns of  $P_{\sharp}$. Then we observe the following matrix invertibility preserving structure: 
 \small
\begin{equation*}
     \Lambda_{n_{c}}=\begin {bmatrix}
\alpha-i\beta & 0 & 0 & 0\\
0 & \alpha+i\beta & 0 & 0 \\
0 & 0 & \gamma+i\delta & 0 \\
0 & 0 & 0 & \gamma-i\delta 
\end{bmatrix},
 P_{\sharp}^{*}MP_{\sharp}=\begin {bmatrix} 0 & a-ib & 0 & 0 \\
a+ib & 0 & 0 & 0 \\
0 & 0 & 0 & c+id\\
0 & 0 & c-id & 0
\end{bmatrix},
\end{equation*}
\normalsize
which clearly shows that $(P_{\sharp}^{*}MP_{\sharp})^{-1}$ exists, and automatically $P_{\sharp}^{*}M S_{\sharp} =0$. On the contrary, if we consider three complex eigenvalues and their corresponding eigenvectors as the columns of $P_{\sharp}$ then we observe
\begin{equation*}
 \Lambda_{n_{c}}=\begin {bmatrix}
\alpha+i\beta & 0 & 0  \\
0 & \alpha-i\beta & 0\\
0 & 0 & \gamma+i\delta \\
\end{bmatrix},
  P_{\sharp}^{*}MP_{\sharp}=\begin {bmatrix}
0 & a-ib & 0  \\
a+ib & 0 & 0\\
0 & 0 & 0 \\
\end{bmatrix}.
\end{equation*}
\normalsize
As we have not included the eigenvector corresponding to the complex conjugate eigenvalue $\gamma-i\delta$ in the column of $P_{\sharp}$, we see that $P_{\sharp}^{*}MP_{\sharp}$ is singular. Since the orphaned eigenvector corresponding to the eigenvalue $\gamma-i\delta$ is now included automatically in first column of $S_{\sharp}$ and the eigenvector corresponding to the eigenvalue $\gamma+i\delta$ is the last column of $P_{\sharp}$, we also have that $P_{\sharp}^{*}MS_{\sharp} \ne 0$.
\end{remark}

\section{Numerical examples}\label{sec:num} In support of the proposed theory, we present numerical results in this section. For nonsymmetric problems, we use the 2D constant direction advection-diffusion problem and a random sparse matrix; for symmetric indefinite problems that represent saddle-point systems, we use the 2D Wilson Dirac system and a mixed Darcy flow problem. To carry out the numerical experiments, we employ the power iteration method \cite{demmel1997applied} to estimate $\rho(E_{TG})$.  For each experiment (i.e., each data point in a plot), the operator $E_{TG}$ is constructed by first explicitly using the $n_c$ smallest generalized eigenvectors to build $P_{\sharp}$ and $R_{\sharp}$, which are then combined with $A$ and $M$ to construct $E_{TG}$ from equation (\ref{ETG_nonsym}). Specifically, we analyze both the error reduction factor and the residual reduction factor for $E_{TG}$ as approximations to the predicted asymptotic convergence rate of $E_{TG}$ from Corollary \ref{cor:optimality_transfer}. The power method is initiated with a randomly generated error vector, denoted by $\mathbf{e}_{k}$. During each iteration, we update the error vector by applying $E_{TG}$, compute the norm ratio of the updated error to the previous error, and normalize the updated error vector. The error reduction factor, defined as $\varrho_{k}=\frac{||\mathbf{e}_{k+1}||}{||\mathbf{e}_{k}||}$, and the residual reduction factor, defined as $\varrho_{A_{k}}=\frac{||A\mathbf{e}_{k+1}||}{||A\mathbf{e}_{k}||}=\frac{||\mathbf{r}_{k+1}||}{||\mathbf{r}_{k}||}$, are recorded for each iteration. Through this iterative process, both $\varrho_{k}$ and $\varrho_{A_{k}}$ converge towards the spectral radius of $E_{TG}$, providing valuable insights into the convergence behavior of the two-grid method. Lastly, for various choices of coarse points $n_{c}$, we compare those  computed factors with the theoretical identity $|1 - \lambda_{n_{c}+1}|$ from Corollary \ref{cor:optimality_transfer}, with the expectation that they match. We use $1500$ power iterations for all nonsymmetric test cases and $500$ power iterations for indefinite problems (Dirac and Darcy) to compute the error and residual reduction factors, which based on our experiments is sufficient iterations to accurately measure the asymptotic rate. 

\subsection{2D constant direction advection-diffusion} \label{sec:adv-diff} First we consider a nonsymmetric example. Here the PDE is
\begin{subequations}
\begin{align}
-\alpha \nabla \cdot \nabla u + \mathbf{b}(x,y) \cdot \nabla u &= f \quad \mbox{for } \Omega = [-1,1]^2, \label{eqn:advdiff1}\\
u &= 0 \quad \mbox{on } \partial \Omega, \label{eqn:advdiff2}
\end{align}
\end{subequations}
where $\alpha$ is the diffusion constant and $\mathbf{b}(x,y)$ is the advection direction.  For this problem, we consider constant non-grid aligned advection $\mathbf{b}(x,y) = [ \sqrt{2/3},\, \sqrt{1/3} ]$. The discretization is first-order upwind discontinuous Galerkin (DG) \cite{DGbook} for advection and the interior-penalty method \cite{arnold1982interior} for the DG discretization of diffusion. PyMFEM \cite{mfem,pymfem} examples 1 (diffusion) and 9 (advection) were used to generate the discretization. We consider the diffusion coefficients $\alpha=10^{-5}$ and $\alpha=10$, representing both advection- and diffusion-dominated problems, respectively, to demonstrate that theory holds for both types of problems. 

\begin{figure}[h!]
     \centering
     \begin{subfigure}[b]{0.49\textwidth}
         \centering
         \includegraphics[width=0.8\textwidth]{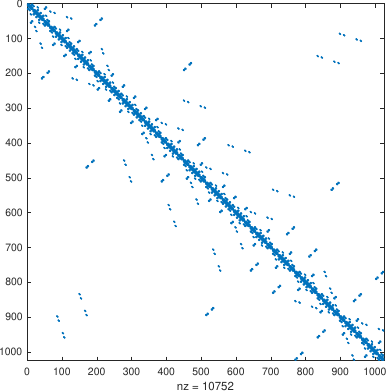}
         \caption{Sparsity pattern of the steady-state nonsymmetric advection diffusion operator $A$ .}
         \label{fig:adv_diff_1}
     \end{subfigure}
     \hfill
     \begin{subfigure}[b]{0.49\textwidth}
         \centering
         \includegraphics[width=\textwidth]{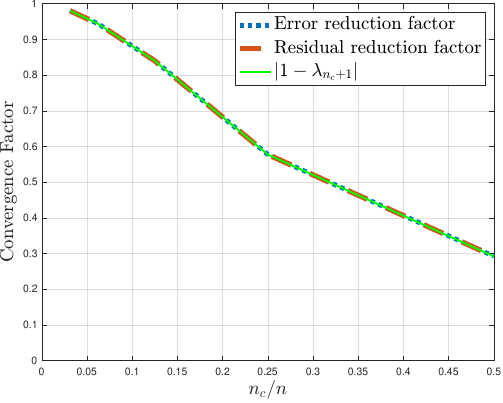}
         \caption{Advection dominant(small diffusion constant $\alpha=10^{-5}$) problem.}
         \label{fig:adv_diff_2}
     \end{subfigure}
     \hfill
     \begin{subfigure}[b]{0.49\textwidth}
         \centering
         \includegraphics[width=\textwidth]{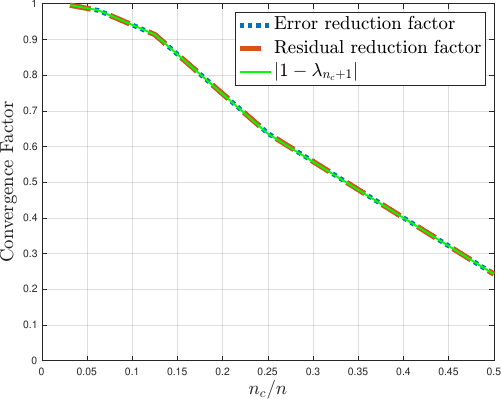}
         \caption{Diffusion dominant (large diffusion constant $\alpha=10$) problem.}
         \label{fig:adv_diff_3}
     \end{subfigure}
     \hfill
     \begin{subfigure}[b]{0.49\textwidth}
         \centering
         \includegraphics[width=\textwidth]{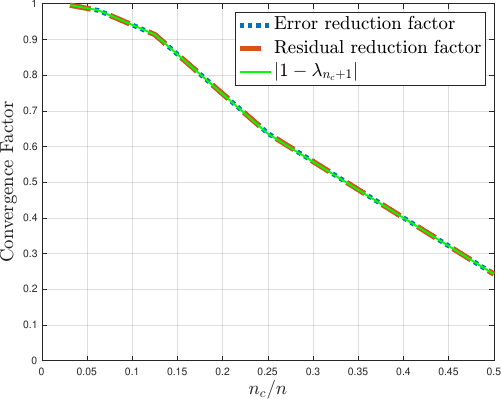}
         \caption{Advection-diffusion mixed (diffusion constant $\alpha=0.1$) problem.}
         \label{fig:adv_diff_nonsym}
     \end{subfigure}
        \label{fig:advection_diffusion}
        \caption{Advection-diffusion problems \eqref{eqn:advdiff1},\eqref{eqn:advdiff2} -  sparsity pattern and verification of the generalized optimal convergence theory using Kaczmarz relaxation where we compare the theoretical prediction with the computed reduction factors. }
\end{figure}

For this example, the discretization matrix $A$ here is nonsymmetric but is definite in a field-of-values sense, with eigenvalues that have nonnegative real part. The sparsity pattern of the nonsymmetric operator $A$ is shown in Figure \ref{fig:adv_diff_1}.
The agreement between theory and the computed reduction factors is shown in Figure \ref{fig:adv_diff_2} for the advection dominant problem and in Figure \ref{fig:adv_diff_3} for the diffusion dominant problem, where it is evident that both reduction factors coincide with the theoretical result $|1 - \lambda_{n_{c}+1}|$, and that the convergence factors decrease with an increase in the number of coarse points $n_{c}$. We also consider the advection-diffusion problem \eqref{eqn:advdiff1} with a diffusion coefficient $\alpha = 0.1$ as an illustrative example to obtain a nonsymmetric operator $A$ representing the advection-diffusion mixed problem. Figure \ref{fig:adv_diff_nonsym} illustrates the agreement between theory and reduction factors for this case. The nonsymmetric operator $A$ in all these examples has dimensions $1024 \times 1024$. For all test cases, we carefully choose the values of $n_{c}$ to be $[32,64,128,256,512]$ in order to ensure that when a complex conjugate eigenvalue pair occurs for the generalized eigenvalue problem of the matrix stencil $(A,M)$, both the complex conjugate right and left eigenvector pairs are placed in $P$ and $R$ respectively. For nonsymmetric matrices, this is not necessary (\Cref{lem:orth} vs. \Cref{cor:hermitian_orthogonality}), but since conjugate pairs will share convergence factor in \Cref{th:optimalPR_nonsymmetric}, it does not make sense to include one and not the other.

\begin{figure}[h!]
     \centering
     \begin{subfigure}[b]{0.49\textwidth}
         \centering
         \includegraphics[width=\textwidth]{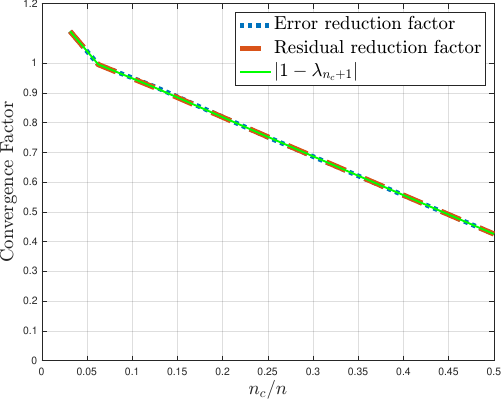}
         \caption{Comparison between theoretical prediction and computed reduction factors using Jacobi relaxation.}
         \label{fig:adv_diff_4}
     \end{subfigure}
     \hfill
     \begin{subfigure}[b]{0.49\textwidth}
         \centering
         \includegraphics[width=\textwidth]{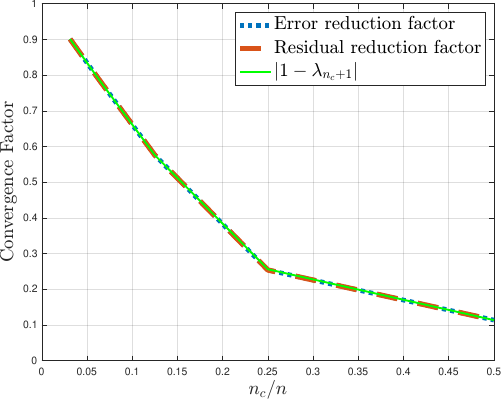}
         \caption{Comparison between theoretical prediction and computed reduction factors using Gauss-Seidel Relaxation.}
         \label{fig:adv_diff_5}
     \end{subfigure}
        \caption{Verification of the generalized optimal convergence theory using Jacobi and Gauss-Seidel method for the advection-diffusion problem \eqref{eqn:advdiff1},\eqref{eqn:advdiff2} with diffusion coefficient $\alpha=0.1$.}
           \label{fig:advection_diffusion_jacobi_gs}
\end{figure}

Additionally, we test our proposed theory identity using Jacobi and Gauss-Seidel relaxation methods. We again select the advection-diffusion problem \eqref{eqn:advdiff1} with diffusion coefficient $\alpha=0.1$ representing a mixed problem in order to show that the suggested theory also applies for these relaxations. Results are shown in Figure \ref{fig:advection_diffusion_jacobi_gs}, where the dimension for nonsymmetric operator $A$ in this example is again $1024\times 1024$. The identity $| 1 - \lambda_{n_{c}+1} |$ is verified for different choices of $n_{c}=32,64,128,256,512$.

Table \ref{table1} presents a comparison of the convergence factors between three different relaxations for this 2D constant advection-diffusion mixed problem \eqref{eqn:advdiff1} with diffusion coefficient $\alpha=0.1$. From Table \ref{table1}, we clearly see that as we increase $n_{c}$, the best convergence factor can be achieved using the Gauss-Seidel relaxation $M$. In fact, asymptotic convergence of  Gauss-Seidel is $27.38\times$ faster than Kaczmarz for aggressive coarsening with $n_c = 32$. For $n_c=64$ and $n_c=128$ Gauss-Seidel is $12.69\times$ and $6.08\times$ faster, respectively, demonstrating the practical difficulties with Kaczmarz relaxation (in addition to added computational cost), despite theoretical guarantees of convergence. 
For Jacobi with ${n_{c}}=32$, we see that the identity $| 1 - \lambda_{n_{c}+1} |>1$ indicates the divergence of the two-grid method. This happens because $\rho (I-M^{-1}A)>1$, that is, the relaxation is divergent here. As we increase $n_{c}$, we eventually obtain a two-grid convergence factor $| 1 - \lambda_{n_{c}+1} |<1$.
\begin{table}[h!]
\begin{center}
\begin{tabular}{ |c|c|c|c|c|c|} 
\hline
$\bm{n_{c}}$ & \textbf{Kaczmarz} & \textbf{Jacobi} & \textbf{GS} & \textbf{Speedup(GS/K)} & \textbf{Speedup(GS/J)} \\ 
\hline
32  & 0.9962  & 1.1074  & 0.9010  & 27.38  & Jacobi Diverges \\ \hline
64  & 0.9816  & 0.9958  & 0.7900  & 12.69 & 56.00 \\ \hline
128 & 0.9132  & 0.9179  & 0.5755  & 6.08  & 6.45 \\ \hline
256 & 0.6356  & 0.7521  & 0.2558  & 3.01  & 4.78 \\ \hline
512 & 0.2421  & 0.4255  & 0.1133  & 1.74  &  2.04\\ \hline
\end{tabular}
\end{center}
\caption{Comparison of asymptotic convergence factors in terms of the spectral radius $\rho(E_{TG}(P,R))$ for three different relaxations $M$ and for various choices of $n_{c}$. Speedups, $\text{ (GS/X)} = \frac{\ln(1/\rho_{\text{GS}})}{\ln(1/\rho_{\text{X}})}$, of Gauss-Seidel (GS) over Kaczmarz (K) and Jacobi (J) are also presented. Specific cases where Jacobi diverges are explicitly mentioned.}
\label{table1}
\end{table}

\subsection{Sparse random matrix} \label{sec:random} 

\begin{figure}[h!]
     \centering
     \begin{subfigure}[b]{0.49\textwidth}
         \centering
         \includegraphics[width=0.8\textwidth]{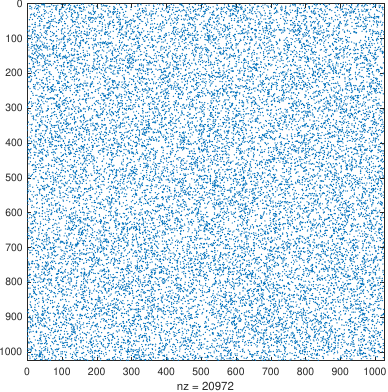}
         \caption{Sparsity pattern of a randomly generated sparse nonsymmetric operator $A$ of size $1024 \times 1024$ with a density of $0.02$.} 
         \label{fig:random_sparse_1}
     \end{subfigure}
     \hfill
     \begin{subfigure}[b]{0.49\textwidth}
         \centering
         \includegraphics[width=\textwidth]{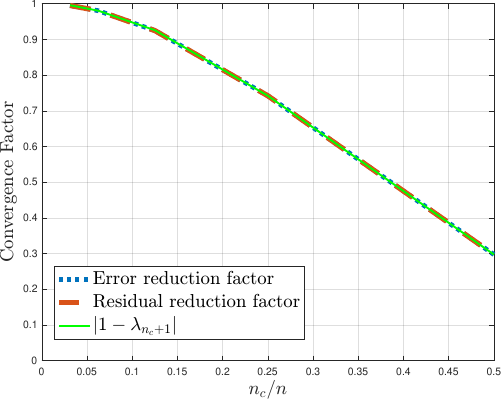}
         \caption{Comparison between theoretical prediction and computed reduction factors using Kaczmarz Relaxation.}
         \label{fig:random_sparse_2}
     \end{subfigure}
        \caption{Verification of the generalized optimal convergence theory using Kaczmarz relaxation method for randomly generated sparse nonsymmetric operator $A$.}
           \label{fig:random_nonsymmetric_experiment}
\end{figure}

To demonstrate that the proposed theory is not restricted to the specific nonsymmetric example presented earlier, we extended our numerical experiments to include random sparse nonsymmetric operators. These matrices were generated using a uniform distribution with a density (proportion of non-zero elements within the matrix) of $0.02$ created in MATLAB using\newline 
\texttt{$\phantom{aaa}$ A=sprand(row,column,density)}\newline 
where \texttt{row,column} denote the matrix size and \texttt{density} is 0.02.
The results confirm that the theoretical predictions and reduction factors remain consistent across varying problem sizes. As an  example, Figure \ref{fig:random_sparse_1} shows the sparsity pattern of a randomly generated sparse nonsymmetric matrix of size $1024 \times 1024$, while Figure \ref{fig:random_sparse_2} demonstrates the agreement between the computed reduction factors and the theoretical identity. These findings further validate the broad applicability of the proposed theory. 

\subsection{2D Wilson Dirac System} Quantum chromodynamics (QCD) on a lattice is a numerical method for calculating quark observables, or elementary particle properties \cite{kogut1983lattice,degrand2006lattice}. In order to approximate the QCD path integral for quark simulations, discrete realizations of the gauge fields must be created. The observables are then computed by averaging over these ensembles of configurations.  For many realizations of the gauge fields, the discretized Dirac equation
\begin{equation}\label{eq:dirac}
  \mathcal{D} \psi = b
\end{equation} must be solved in each of these phases of a lattice QCD calculation. The Bootstrap AMG (BAMG) approach has been used \cite{brannick2014bootstrap} previously to successfully solve such a problem. We consider Wilson's discretization~\cite{wilson} of the Dirac
equation so that $\mathcal{D} = \mathcal{D}_{0} +mI$, where the non-Hermitian massless Wilson matrix is represented by $\mathcal{D}_{0}$, and the quarks' mass is indicated by the shift $m$ \cite{sharpe2006discretization, splittorff2011wilson}. We refer to $\mathcal{D}$ as just the Wilson matrix, which suggests that it is related to a specific shift $m$. Eigenvalues of the Wilson matrix $\mathcal{D}$ are very close to zero, which leads to the highly ill-conditioned non-Hermitian matrix.

\begin{figure}[h!]
     \centering
     \begin{subfigure}[b]{0.49\textwidth}
         \centering
\includegraphics[width=0.8\textwidth]{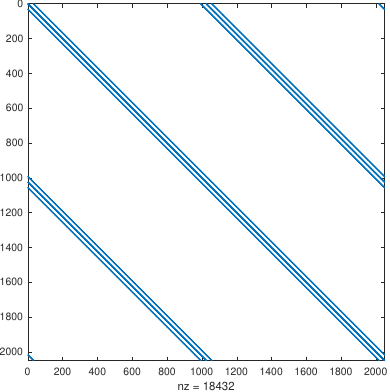}
         \caption{Sparsity pattern of the Hermitian indefinite Wilson matrix $Z$, a representative example for a saddle-point system.}
         \label{fig:dirac_1}
     \end{subfigure}
     \hfill
     \begin{subfigure}[b]{0.49\textwidth}
         \centering
         \includegraphics[width=\textwidth]{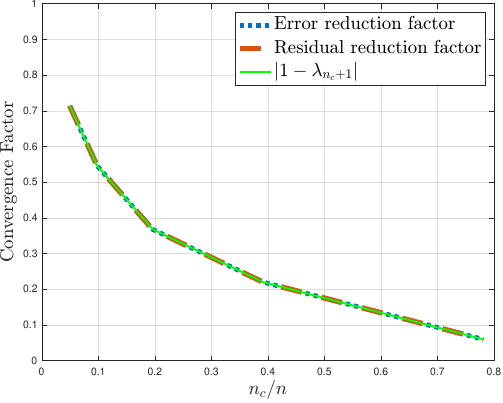}
         \caption{Comparison between theoretical prediction and computed two-grid reduction factors using symmetrized Kaczmarz relaxation.  }
         \label{fig:dirac_2}
     \end{subfigure}
     \hfill
     \begin{subfigure}[b]{0.49\textwidth}
         \centering
         \includegraphics[width=\textwidth]{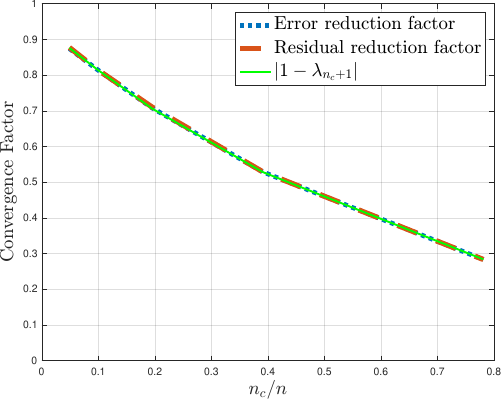}
         \caption{Comparison between theoretical prediction and computed two-grid reduction factors using symmetrized Jacobi relaxation.  }
         \label{fig:dirac_3}
         \end{subfigure}
         \hfill
     \begin{subfigure}[b]{0.49\textwidth}
         \centering
         \includegraphics[width=\textwidth]{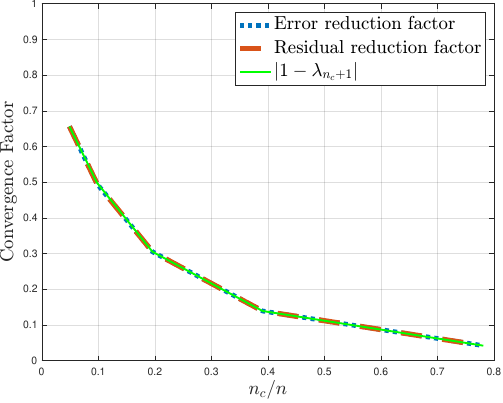}
         \caption{Comparison between theoretical prediction and computed two-grid reduction factors using symmetrized Gauss-Seidel relaxation.  }
         \label{fig:dirac_4}
         \end{subfigure}
    \label{fig:dirac_equation}
    \caption{For the Wilson discretization of the Dirac equation problem \eqref{eq:dirac} - sparsity pattern and numerical verification of the generalized optimal convergence theory.}
\end{figure}

The $\gamma_{5}$-symmetry, resulting from the Clifford algebra \cite{ablamowicz2004applications} structure of the $\gamma$-matrices, is inherited by the Wilson matrix, which is a discretization of the Dirac equation \eqref{eq:dirac}. In particular, $\gamma_{5}:=i\gamma_{1}\gamma_{2}=\operatorname{diag}(-1,1)$ symmetrizes the Dirac equation where
\begin{equation*}
    \gamma_{1}=\begin{pmatrix}
        0 & 1 \\
        1 & 0
    \end{pmatrix} \,\,\,\text{and}\,\,\,
    \gamma_{2}=\begin{pmatrix}
        0 & i \\
        -i & 0
    \end{pmatrix}.
\end{equation*} As a result, we define $\Gamma_{5} = I_{N^{2}}\otimes \gamma_{5}$ in the discrete setting where $N$ is the number of grid points in a given space-time dimension. Scaling the non-Hermitian Wilson matrix $\mathcal{D}$ by $\Gamma_{5}$, we obtain the Hermitian and indefinite form of the Wilson matrix $Z=\Gamma_{5}\mathcal{D}$, where $\Gamma_{5}$ is a unitary matrix \cite{babich2010adaptive}. To numerically verify the optimal interpolation theory, we take  $N=32$, which results in $2N^{2}\times 2N^{2}$ degrees of freedom for both the non-Hermitian Wilson matrix $\mathcal{D}$ and its Hermitian indefinite variant $Z$. The theoretical verification for the Hermitian indefinite Wilson matrix $Z$, where we utilize a symmetrized form of Kaczmarz relaxation $\widetilde{M}$ using \eqref{sym_smoother1}, is shown in Figure \textcolor{red}{\ref{fig:dirac_2}}. We also verify the theory using symmetrized Jacobi and Gauss-Seidel Relaxation $\widetilde{M}$ in Figure \ref{fig:dirac_3} and Figure \ref{fig:dirac_4}, respectively.  

Table \ref{table2} presents the asymptotic convergence factors, measured by the spectral radius $\rho(E_{TG}(P_\sharp))$, for the Hermitian indefinite Wilson matrix $Z$ under three symmetrized relaxation schemes: Kaczmarz, Jacobi, and Gauss-Seidel. It also includes speedups of Gauss-Seidel over Kaczmarz and Jacobi, quantified by the ratio of logarithmic convergence factors. As $n_{c}$ increases, the convergence factors improve for all methods, indicating faster convergence for larger $n_{c}$. Gauss-Seidel consistently outperforms Kaczmarz and Jacobi, achieving the highest speedup over Jacobi, particularly for smaller $n_{c}$.
\begin{table}[h!]
\begin{center}
\begin{tabular}{ |c|c|c|c|c|c|} 
\hline
$\bm{n_{c}}$ & \textbf{Kaczmarz} & \textbf{Jacobi} & \textbf{GS} & \textbf{Speedup(GS/K)} & \textbf{Speedup(GS/J)} \\ 
\hline
100  & 0.7141  & 0.8760  & 0.6560  & 1.25  & 3.18 \\ \hline
200  & 0.5461  & 0.8162  & 0.4980  & 1.15 & 3.43 \\ \hline
400 & 0.3682  & 0.7061  & 0.3048  & 1.19  & 3.41 \\ \hline
800 & 0.2199  & 0.5284  & 0.1388  & 1.30  & 3.09 \\ \hline
1600 & 0.0593  & 0.2836  & 0.0416  & 1.12  &  2.52\\ \hline
\end{tabular}
\end{center}
\caption{Comparison of asymptotic convergence factors for Hermitian indefinite Wilson matrix $Z$ in terms of the spectral radius $\rho(E_{TG}(P_\sharp))$ for three different symmetrized relaxations $\widetilde{M}$ and for various choices of $n_{c}$. Speedups, $\text{ (GS/X)} = \frac{\ln(1/\rho_{\text{GS}})}{\ln(1/\rho_{\text{X}})}$, of Gauss-Seidel (GS) over Kaczmarz (K) and Jacobi (J) are also presented.}
\label{table2}
\end{table}
\subsection{2D mixed Darcy problem} The behavior of fluids in porous media is described by the Darcy equations. These comprise a set of PDEs that come from Darcy's law and conservation of mass. A common method for discretizing the Darcy equations and obtaining an approximate solution is the use of mixed finite element methods. This falls into the category of numerical techniques that rely on a variational formulation of the problem, where the variables are pressure and velocity. In this example, we  consider a simple 2D mixed Darcy problem, 

\begin{subequations}
\begin{align}
k\,\,\mathbf{u} + \nabla p = f \label{eqn:darcy1}\quad \mbox{in } \Omega\\
 - \nabla \cdot \mathbf{u} = g  \quad \mbox{in } \Omega ,\label{eqn:darcy2}
\end{align}
\end{subequations}
where $\Omega\subset\mathbb{R}^d$ is a bounded connected flow domain, i.e., a porous media saturated with a fluid, with $-p= \text{``given pressure"}$ as its natural boundary condition in the boundary $\partial\Omega$. In \eqref{eqn:darcy1} and \eqref{eqn:darcy2}, $\mathbf{u}:\Omega \to \mathbb{R}^d$ is the fluid velocity, $p:\Omega \to \mathbb{R}$ is the fluid pressure, $f,g \in L^{2}(\Omega)$ represents potential source (or sinks) in the medium and $k$ is the permeability constant of the medium. For the  fluid velocity $\mathbf{u}$, the first equation \eqref{eqn:darcy1} is Darcy's law, and the second equation \eqref{eqn:darcy2}, often known as the continuity equation, is derived from the law of conservation of mass.
The discretization is $H(\operatorname{div},\Omega)$ Raviart-Thomas
finite elements \cite{raviart_thomas} for velocity $\mathbf{u}$ and $L^{2}(\Omega)$ piecewise discontinuous
 polynomials for pressure $p$. For this example, we use constant coefficient $k=1.0$, and a given exact solution $\mathbf{u}=(u_{x},u_{y})=(- e^{x}\sin{y},- e^{x}\cos{y}), p=e^{x}\sin{y}$ to compute the corresponding right hand sides $f$ and $g$.  Using PyMFEM \cite{mfem,pymfem} example 5 (Darcy Problem), we generate the discretizations. 
 \begin{figure}[h!]
     \centering
     \begin{subfigure}[b]{0.35\textwidth}
         \centering
\includegraphics[width=\textwidth]{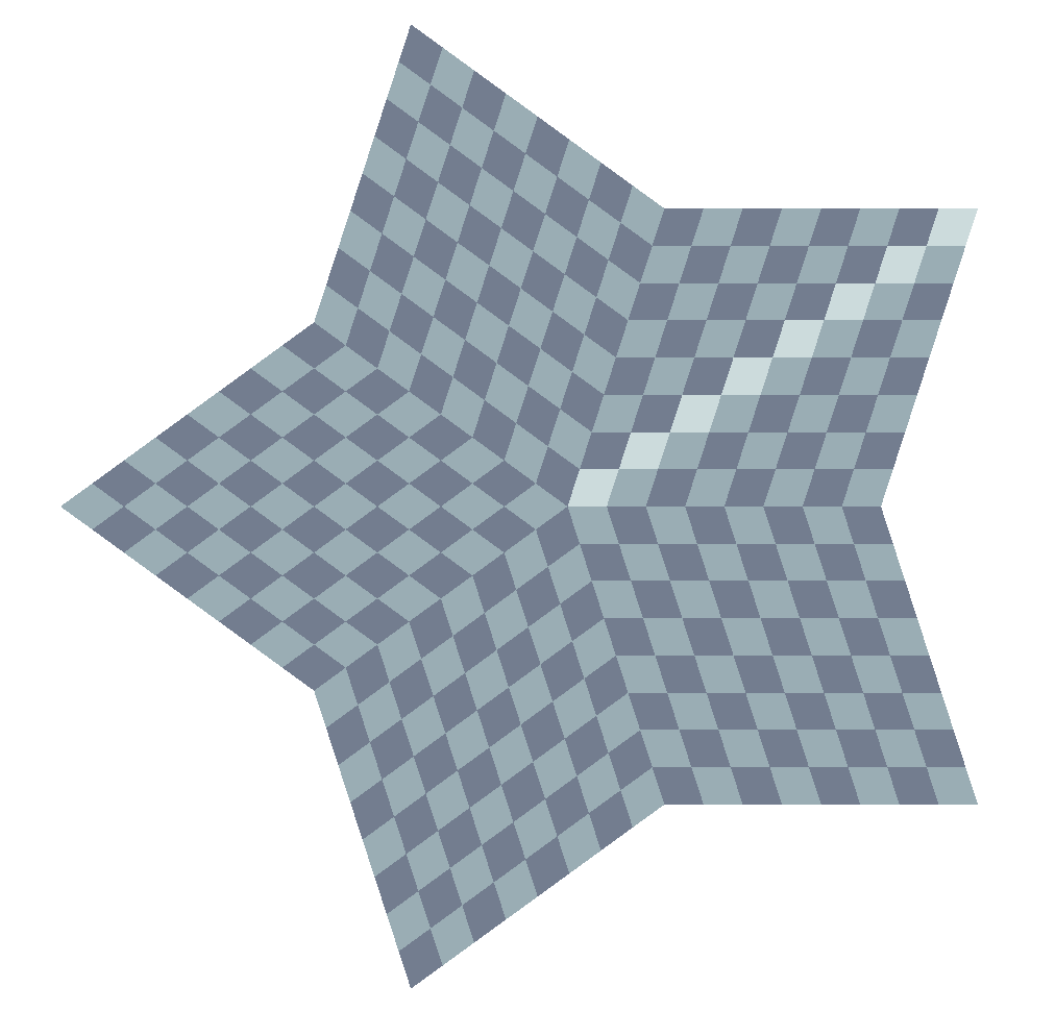}
         \caption{Star mesh used in finite element discretization}
         \label{fig:darcy_1}
     \end{subfigure}
     \hfill
     \begin{subfigure}[b]{0.49\textwidth}
         \centering
         \includegraphics[width=0.8\textwidth]{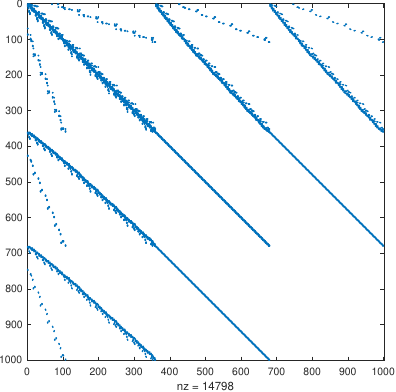}
         \caption{Sparsity pattern of the Darcy operator $\mathcal{A}$  }
         \label{fig:darcy_2}
     \end{subfigure}
     \hfill
     \begin{subfigure}[b]{0.49\textwidth}
         \centering
         \includegraphics[width=\textwidth]{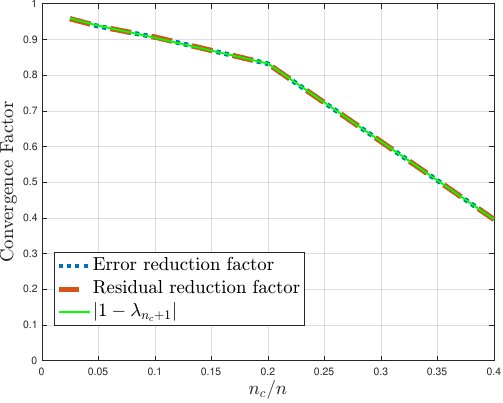}
         \caption{Comparison between theoretical prediction and computed reduction factors using symmetrized Kaczmarz relaxation $\widetilde{M}$.}
         \label{fig:darcy_3}
     \end{subfigure}
    \label{fig:darcy_equation}
    \caption{For the 2D mixed Darcy problem \ref{eqn:darcy1}, \ref{eqn:darcy2} - Initial mesh, sparsity pattern and numerical verification of the generalized optimal convergence theory.}
\end{figure}
Mixed finite element discretizations of the Darcy problem, \eqref{eqn:darcy1}-\eqref{eqn:darcy2}, lead to linear systems of equations in saddle-point form. We denote this symmetric indefinite Darcy operator by 
 \begin{equation}\label{darcy_assemble}
     \mathcal{A}=\begin{pmatrix}
         K & B^{*}\\
         B  &  0
     \end{pmatrix},
 \end{equation}
 where \begin{align*}
     K&=\displaystyle\int_\Omega k\, u_{h} \cdot v_{h} \,d\Omega,  \quad  u_{h}, v_{h} \in V_{h};\\
     B&=-\displaystyle\int_\Omega (\nabla\cdot u_{h}) q_{h} \,d\Omega,\quad   u_{h} \in V_{h}, q_{h} \in Q_{h};
 \end{align*}
 $V$ and $Q$ are two Hilbert spaces with their corresponding finite dimensional sub-spaces,  $V^{h}$ and $Q^{h}$, respectively. Using the Darcy operator \eqref{darcy_assemble}, we verify the generalized optimal interpolation-based AMG convergence theory. Figure \ref{fig:darcy_1} depicts the star-shaped mesh used in the finite element discretization of the 2D mixed Darcy problem, which serves as the foundation for assembling the operator. The sparsity pattern of the resulting Darcy operator $\mathcal{A}$ is shown in Figure \ref{fig:darcy_2}, illustrating the structure of the discretized system. Finally, Figure \ref{fig:darcy_3} compares the theoretical prediction and computed reduction rates using the symmetrized Kaczmarz relaxation $\widetilde{M}$, demonstrating that the convergence identity holds for this problem as well. Note that, for this example, the direct implementation of Jacobi or Gauss-Seidel as a relaxation operator results in a singular matrix $M$. Consequently, the theory was verified exclusively using the symmetrized Kaczmarz operator $\widetilde{M}$.
 \begin{figure}[h!]
     \centering
     \begin{subfigure}[b]{0.32\textwidth}
         \centering
\includegraphics[width=\textwidth]{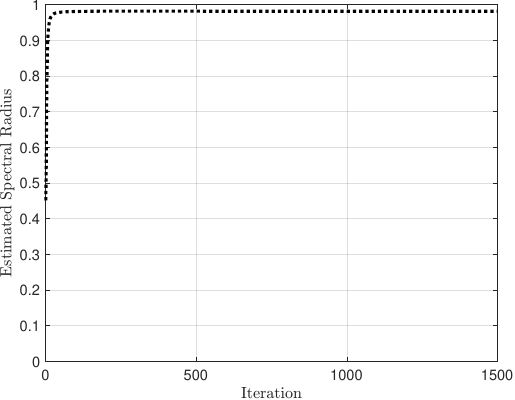}
         \caption{Advection-diffusion}
         \label{fig:power_1}
     \end{subfigure}
     \hfill
     \begin{subfigure}[b]{0.32\textwidth}
         \centering
         \includegraphics[width=\textwidth]{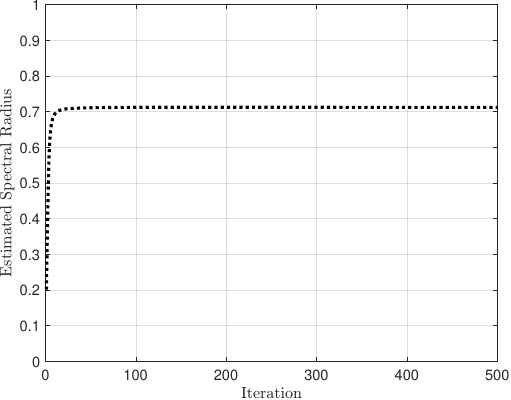}
         \caption{Dirac }
         \label{fig:power_2}
     \end{subfigure}
     \hfill
     \begin{subfigure}[b]{0.32\textwidth}
         \centering
         \includegraphics[width=\textwidth]{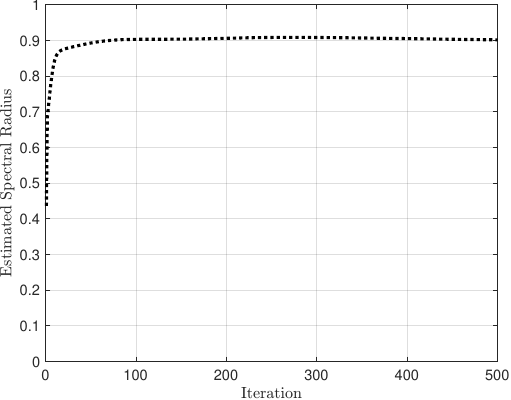}
         \caption{Darcy}
         \label{fig:power_3}
     \end{subfigure}
    \caption{Asymptotic convergence of the power iteration for three different examples is shown. Figure \ref{fig:power_1} illustrates the advection-diffusion problem defined by equations \eqref{eqn:advdiff1} and \eqref{eqn:advdiff2} with a diffusion constant $\alpha = 0.1$, where $n = 1024$ and $n_{c} = 64$. Figure \ref{fig:power_2} depicts the Dirac problem as described by equation \eqref{eq:dirac}, with $n = 2048$ and $n_{c} = 101$. Figure \ref{fig:power_3} presents the Darcy equation problem given by equations \eqref{eqn:darcy1} and \eqref{eqn:darcy2}, with $n = 2048$ and $n_{c} = 100$. The advection-diffusion problem in Figure \ref{fig:power_1} employs Kaczmarz relaxation $M$, while Figures \ref{fig:power_2} and \ref{fig:power_3} utilize symmetrized Kaczmarz relaxation $\widetilde{M}$.}
    \label{fig:power_iteration}
\end{figure}
\begin{remark}[On asymptotic nature of spectral radius]
It is well-known that spectral radius is only a guaranteed asymptotic measure of convergence for nonsymmetric problems, whereas we are interested in convergence rates in $\mathcal{O}(1)$ iterations. For example, one could have a nilpotent matrix with all zero eigenvalues that exhibits exponential growth in solution magnitude for some finite number of iterations before decaying rapidly as nilpotency sets in. Figure \ref{fig:power_iteration} plots the evolution of our power iteration measuring error and residual decay for each of the three problems discussed previously. Remember, a smaller value on the $y$-axis indicates faster convergence.  Notably, we see no  misleading behavior from the power method, i.e., it is monotonically increasing.  In other words, the asymptotic behavior of $E_{TG}$ (during later power iterations) provides a good upper bound on the convergence behavior of $E_{TG}$ during earlier iterations. Thus for these problems, the spectral radius provides a robust and practical estimate of convergence.
\end{remark}

\section{Conclusion}\label{sec:conc} In this paper, we  present a generalization of the two grid optimal AMG convergence theory that can be applied to general nonsymmetric problems. Moreover, we show that this theory result can be naturally extended to symmetric indefinite problems. In order to support the derived convergence identity based on optimal interpolation and restriction operators, we offer a straightforward proof for the proposed theory which relies on a specific measure of spectral radius as well as numerical evidence. This new conclusion gives us an identity that holds for any matrix and specifies the asymptotic convergence rate of a two-grid method. This is the first of its kind to provide a theoretical convergence rate of a two-grid approach for general matrices, to our knowledge. Moreover by establishing a theoretical framework that rigorously links the transfer operators and relaxation schemes, this theory will allow for more detailed investigations and analyses of existing and future nonsymmetric multigrid solvers, thus aiding in the development of more robust nonsymmetric solvers. Towards this goal, future work will be focused on creating an AMG solver for nonsymmetric and indefinite saddle-point systems using a generalized reduction approach, where we will take into account the proposed theory when designing this solver, which will provide a target convergence for a particular coarsening scheme. This ability to evaluate observed convergence rates for practical AMG methods, in comparison to this theory, is a key practical contribution of this work.

\bibliographystyle{siamplain}
\bibliography{references}
\end{document}